\documentclass[11pt]{article}

\usepackage[margin=1in]{geometry}
\usepackage{amsmath}
\usepackage{amssymb}
\usepackage{graphicx} 
\usepackage{mathtools}
\usepackage{amsthm}
\usepackage[colorlinks=true,linkcolor=blue,citecolor=blue]{hyperref}
\usepackage{multirow}
\theoremstyle{plain}
\newtheorem{theorem}{Theorem}
\newtheorem{lemma}{Lemma}

\newtheorem{definition}{Definition}
\usepackage{natbib}

\newtheorem{remark}{Remark}

\newcommand{\R}{\mathbb{R}}
\newcommand{\Q}{\mathbb{Q}}
\newcommand{\1}{\mathbf{1}}

\title{Polynomial-Time Solvability of the Multi-Dimensional Slope Change Problem}
\author{Apoorva Narula, Santanu S. Dey\\
H. Milton Stewart School of Industrial and Systems Engineering (ISyE),\\
Georgia Institute of Technology, USA\\
anarula34@gatech.edu, santanu.dey@isye.gatech.edu}
\date{}

\begin{document}

\maketitle

\begin{abstract}
Slope change detection seeks to identify changes in the trends of time-ordered signals. We study the problem of fitting a multidimensional time series by continuous piecewise-linear functions with a common, continuously varying breakpoint and unrestricted slopes. We consider both the dense setting, in which every dimension may change slope, and the sparse setting, in which at most \(k\) dimensions may change slope. Although the resulting optimization problems are nonconvex mixed integer nonlinear programs and may have irrational optimal solutions even for rational input data, we prove that all six variants, corresponding to the \(\ell_1\), \(\ell_2\), and \(\ell_\infty\) losses in the dense and sparse settings, are exactly solvable in polynomial time. In particular, optimal solutions admit polynomial-size algebraic representations and can be computed with polynomial bit complexity. For the least-squares loss, we obtain especially efficient algorithms requiring \(O(JT)\) operations in the dense setting and \(O(J^2T\log J)\) operations in the sparse setting, where \(T\) is the number of time points and \(J\) is the number of dimensions. We implement these least-squares algorithms and evaluate them on simulated data and on three real-world multidimensional time series, namely, regional surface temperatures (ERA5), hand-tracking data from a mixed-reality headset (HoloAssist), and turbofan engine sensor data (NASA C-MAPSS). The dense algorithm solves instances with up to \(J=10^5\) and \(T=100\), that is, \(10^7\) observations, in less than \(0.1\) seconds on a laptop, while the sparse algorithm solves instances with about \(10^3\) dimensions and up to \(T=240\) in under a second. Compared with a mixed-integer nonlinear programming (MINLP) formulation solved by Gurobi on the same instances, our algorithms are several orders of magnitude faster (approximately $8.9 \times 10^3$ to $3.4 \times 10^6$ times) on those instances where Gurobi terminates within its default optimality tolerance of $0.01\%$ optimality gap. On larger instances, Gurobi frequently reaches a 10,000-second time limit with substantial optimality gaps.

\end{abstract}


\section{Introduction}

Slope change detection seeks the location in a time-ordered sequence, called a \emph{breakpoint}, at which the rate of change of an underlying signal changes. Unlike mean-shift detection, which seeks abrupt changes in the signal level, slope change detection considers changes in trend while preserving continuity. Slope-change models have been used in a variety of applications, including surface-temperature trends in climatology \cite{A2025Change,rahmstorf2017global,beaulieu2024recent,triacca2026recent}, rainfall trends \cite{de2010breakpoints}, changes in health outcomes measured repeatedly over time \cite{brilleman2017bayesian}, and sensor-based monitoring \cite{cao2018multi}. In many applications, multiple signals are observed at each observation time, and a common event, such as a regime shift or a system fault, may change their trends simultaneously. This motivates fitting a slope change that is shared across the dimensions of a multidimensional time series.

\subsection{Problem formulation}

Let $T\geq4$ be an integer, let $x_1<x_2<\cdots<x_T$ be rational time coordinates, and let $y_i=(y_i^1,\ldots,y_i^J)\in\Q^J$ be the observation at time $x_i$, so that $y_i^j$ is the observation of dimension $j\in[J]$ at time $x_i$. For a breakpoint $z\in\R$, every continuous function that is affine on $(-\infty,z]$ and on $[z,\infty)$ can be written uniquely as
$$x\mapsto\alpha+\beta x+\delta(x-z)_+,$$
where $(t)_+:=\max\{t,0\}$. We call $z$ the breakpoint of this function. This function has slope $\beta$ and intercept $\alpha$ before the breakpoint, and slope $\beta+\delta$ and intercept $\alpha-\delta z$ after it. Thus $\delta$ is the change in slope, and continuity at $z$ is built into the representation. We fit one such function to each dimension, with a breakpoint common to all dimensions. Given $z\in\R$ and $\alpha,\beta,\delta\in\R^J$, the residual of dimension $j$ at time $x_i$ is
\begin{equation}\label{eq:residual}
r_i^j=y_i^j-\alpha_j-\beta_j x_i-\delta_j(x_i-z)_+,\qquad i\in[T],\ j\in[J],
\end{equation}
and $r\in\R^{TJ}$ denotes the vector of these $TJ$ residuals. For $p\in\{1,2,\infty\}$, define the loss
$$L_p(r):=\begin{cases}\sum_{i,j}|r_i^j|, & p=1,\\ \sum_{i,j}(r_i^j)^2, & p=2,\\ \max_{i,j}|r_i^j|, & p=\infty.\end{cases}$$
Given data $(x_i,y_i)_{i\in[T]}$ with $x_i\in\mathbb{Q}$ and $y_i\in\mathbb{Q}^J$, and $k\in\{0,1,\ldots,J\}$, the \emph{slope change problem} is
\begin{equation}\label{eq:slope-change}
\textup{OPT}_{p,k}:=\min_{z,\alpha,\beta,\delta,r}\big\{L_p(r):\ r \text{ satisfies \eqref{eq:residual}},\ \|\delta\|_0\leq k,\ z\in[x_2,x_{T-1}]\big\}.
\end{equation}
The constraint $\|\delta\|_0\leq k$ allows at most $k$ dimensions to change slope, and each remaining dimension is fitted by a single line. We refer to the case $k=J$ as the \emph{dense} problem and to the case $k<J$ as the \emph{sparse} problem. The constraint $z\in[x_2,x_{T-1}]$ ensures that at least two observations lie on each side of the breakpoint, since a change in slope cannot be inferred from a segment containing a single observation.

Three features make this problem nontrivial. First, the breakpoint $z$ is a continuous decision variable, and the terms $\delta_j(x_i-z)_+$ introduce nonconvexity. Second, the breakpoint is shared across all dimensions, coupling their fitting problems. Third, the sparse problem introduces the additional combinatorial problem of selecting which dimensions may change slope. Although the problem can be formulated as a mixed-integer nonlinear program \cite{narula2026changepointdetectionpiecewiselinear}, we show that its special structure permits polynomial-time algorithms.

\subsection{Related literature}

\emph{Single-dimension time series.} Most methods for a single time series allow several breakpoints but restrict the breakpoints to the observation times or a prespecified candidate set. Without continuity, dynamic programming finds exact solutions \cite{jackson2005algorithm,bai2003computation,killick2012optimal}. With continuity, \cite{fearnhead2019detecting} use functional dynamic programming with an $\ell_0$ penalty on the number of slope changes, and \cite{fearnhead2024cpop} allow unevenly spaced observations and more flexible candidate sets. The convex alternative, $\ell_1$ trend filtering \cite{kim2009trend,tibshirani2014adaptive}, penalizes the second differences of the fitted values, so its breakpoints arise implicitly at observation times. Fewer methods treat breakpoint locations as continuous. For one breakpoint under least squares, the setting closest to ours, \cite{hudson1966fitting} obtains the exact fit by examining each interval between consecutive observations, with statistical properties studied in \cite{hinkley1969inference,feder1975asymptotic} and a grid-search alternative proposed by \cite{lerman1980fitting}. For several breakpoints, iterative linearization \cite{muggeo2003estimating} guarantees only local optimality, while mixed-integer programming provides exact formulations \cite{toriello2012fitting,rebennack2020,goldberg2021,narula2026changepointdetectionpiecewiselinear}.

\paragraph{Multiple-dimension time series with shared breakpoints.} Here, the breakpoints are common to all dimensions. \cite{kim2011estimating} estimates a break in the intercept, slope, or both of a linear trend, common to all series in a large panel, by minimizing the total sum of squared residuals over all break dates, and \cite{kim2014common} estimates the break date jointly with latent common factors in the errors. For shared changes in the mean, the group fused lasso \cite{bleakley2011group} gives a convex formulation. \cite{cao2018multi} detects, online, a common change point after which a subset of the dimensions moves from a constant mean to a linear drift, with sparsity modeled through the probability that each dimension is affected. Sparse changes in the mean of high-dimensional time series are studied by \cite{wang2018high} and \cite{enikeeva2019high}. These works differ from our setting in their restrictions on breakpoint locations, assumptions on the form of the slope change, or statistical detection objectives, and none establishes polynomial-time global optimization for the formulation \eqref{eq:slope-change}. \cite{narula2026changepointdetectionpiecewiselinear} extends their mixed-integer formulations to shared breakpoints. These formulations are exact, but they are not known to be solvable in polynomial time.

\paragraph{Parametric linear programming.} \cite{boveroux2026} optimizes the value of $\min\{c^\top x:(A_0+\lambda A_1)x\leq b,\ x\geq0\}$ over a scalar $\lambda$ in an interval, extending sensitivity analysis \cite{freund1985postoptimal}. They prove that both minimizing and maximizing this value are NP-hard, and identify easier cases where $A_1$ has entries of one sign or a single nonzero row or column. For $p\in\{1,\infty\}$ and fixed $i^*$, our dense problem is the minimization version with $\lambda=z$. However, because of the absolute values, $A_1$ always has entries of both signs, and for $J\geq2$ it has $J$ nonzero columns, so none of these cases applies. We nevertheless solve the dense and sparse problems exactly in polynomial time, even though the optimal $z$ may be irrational and the sparse problem is not an LP for fixed $i^*$.


\paragraph{Complexity and connection to ReLU networks.} Closely related function recovery problems are NP-hard \cite{weinmann2015iterative}, and \cite{fearnhead2019detecting} notes that it is unclear whether detecting multiple slope changes is NP-hard. We study the setting with a single breakpoint that is shared across dimensions and located continuously. This setting is also the global training of a ReLU network whose input is time, with a single hidden unit $(x-z)_+$ shared by $J$ outputs with affine skip connections, where the breakpoint $z$ is restricted to $[x_2,x_{T-1}]$. Training ReLU networks is NP-hard for other architectures \cite{dey2020approximation,boob2022complexity,froese2023training}, and optimal weights can be irrational \cite{abrahamsen2021training,bertschinger2023training}. These results do not imply hardness for ours. Exact polynomial-time training is known for one hidden layer with a single output and fixed input dimension \cite{arora2018understanding}, but this does not cover multiple outputs.

\section{Main complexity results}
\label{sec:main}

Throughout, $\langle\mathcal{D}\rangle$ denotes the binary encoding length of the data $\mathcal{D}:=(x_i,y_i)_{i\in[T]}$. The existence and representation of rational solutions are important issues in mixed-integer nonlinear programming. For example, mixed-integer quadratic programming admits polynomial-size certificates of feasibility \cite{delpia2017mixed}, whereas systems of cubic polynomial inequalities with rational data may have only irrational feasible solutions, and deciding the existence of rational solutions is NP-hard \cite{bienstock2023complexity}. Problem \eqref{eq:slope-change} exhibits both phenomena, that is, optimal breakpoints can be irrational, yet we show that optimal solutions admit polynomial-size algebraic representations and can be computed exactly in polynomial time. To see that optimal breakpoints can be irrational, let $T=6$, $J=2$, $x_i=i-1$ for $i\in[6]$, $y^1=(-2,2,-2,-2,5,0)$, and $y^2=(-9,4,8,-6,-4,0)$. Then
$$\textup{OPT}_{2,2}=\frac{74531-\sqrt{477755737}}{420}\approx125.4128,$$
attained at $z^\star=(-58563+5\sqrt{477755737})/49162\approx1.0318$. Since $477755737$ is not a perfect square, both $z^\star$ and $\textup{OPT}_{2,2}$ are irrational. This irrationality is caused by sharing the breakpoint across dimensions. For $J=1$ and $p=2$ a rational optimal breakpoint always exists (see Section~\ref{sec:nonsparse}). We therefore specify what an algorithm returns.

\begin{definition}\label{def:solve}
An algorithm \emph{solves} \eqref{eq:slope-change} if it returns an optimal breakpoint $z^\star$, represented as a real algebraic number by a univariate polynomial with rational coefficients together with an interval with rational endpoints that contains $z^\star$ and no other root of the polynomial, and returns $\alpha^\star,\beta^\star,\delta^\star$ as explicit rational functions of $z^\star$ with rational coefficients.
\end{definition}

As we shall see, for $p\in\{2,\infty\}$, our algorithms return an optimal breakpoint $z^\star$ that is algebraic of degree at most two over $\Q$. For $p=1$, its degree is at most $\max\{2,2J-2\}$ in the dense problem and at most $\max\{2,2k-2\}$ in the sparse problem. In all cases, the optimal fitting coefficients can be chosen in $\Q(z^\star)$, the smallest subfield of $\R$ containing $\Q$ and $z^\star$, which consists of the values at $z^\star$ of rational functions with rational coefficients. 

The encoding length of an algebraic number represented as in Definition~\ref{def:solve} is the total binary encoding length of the coefficients of the polynomial and the endpoints of the interval. All rational and algebraic numbers arising in our algorithms admit representations of encoding length polynomial in $T$, $J$, and $\langle\mathcal{D}\rangle$. Standard algorithms for real algebraic numbers perform exact arithmetic, root isolation, sign evaluation, and comparison on these representations in polynomial time \cite{basu2006algorithms,sagraloff2016computing}. Consequently, the optimal breakpoint, objective value, and fitting coefficients can also be approximated to within $2^{-t}$ in time polynomial in $T$, $J$, $\langle\mathcal{D}\rangle$, and $t$ \cite{BrentZimmermann2010,sagraloff2016computing}.

\subsection{Poly-time complexity}

\begin{theorem}\label{thm:main}
For every $p\in\{1,2,\infty\}$ and $k\in\{0,1,\ldots,J\}$, problem \eqref{eq:slope-change} can be solved, in the sense of Definition~\ref{def:solve}, in time polynomial in $T$, $J$, and $\langle\mathcal{D}\rangle$. More precisely, the algorithms perform the arithmetic operations counted in Table~\ref{tab:complexity}. All numbers involved are rational or algebraic, with encoding length polynomial in $T$, $J$, and $\langle\mathcal{D}\rangle$. For $p\in\{2,\infty\}$, these operations include, and the table counts among them, computations of roots of quadratic polynomials and exact comparisons of quadratic algebraic numbers. For $p=1$, the algorithm additionally performs root isolations and algebraic-number comparisons, counted in Table~\ref{tab:complexity}. All these algebraic operations have polynomial bit complexity.
\end{theorem}








\begin{table}[htbp]
\centering
\caption{Number of operations of the algorithms. Quadratic root computations and comparisons of quadratic algebraic numbers count as single operations. For $p=1$, the table also lists the root isolations of the polynomials of degree at most $2J-2$ or $2k-2$, which are needed only if $J\geq2$ or $k\geq2$, and the comparisons of the resulting algebraic numbers. All these operations have polynomial bit complexity (Theorem~\ref{thm:main}).}
\label{tab:complexity}
\scriptsize
\renewcommand{\arraystretch}{1.3}
\begin{tabular}{c|c|c}
\hline
& \textbf{Dense ($k=J$)} & \textbf{Sparse ($k<J$)}\\
\hline
$p=2$ & $O(JT)$ & $O(J^2T\log J)$\\
\hline
$p=\infty$ & $O(J^3T^9\log(JT))$ & $O(JT+J\log J+k^3T^9\log(kT))$\\
\hline
$p=1$ & \begin{tabular}{@{}c@{}}$O(J^3T^9+JT^9\log(JT))$\\ arithmetic operations;\\ $O(JT^9)$ root isolations,\\ degree $\leq 2J-2$ ($J\geq2$);\\ $O(J^2T^9)$ algebraic comparisons\end{tabular} & \begin{tabular}{@{}c@{}}$O((J+k^2)J^3T^9+JT^9\log(JT))$\\ arithmetic operations;\\ $O(J^3T^9)$ root isolations,\\ degree $\leq 2k-2$ ($k\geq2$);\\ $O(kJ^3T^9)$ algebraic comparisons\end{tabular}\\
\hline
\end{tabular}
\end{table}

When $k=0$, each dimension is fitted by a single line, and the problem is solved directly. For $p=2$, the dense algorithm is linear in the number $JT$ of observations, which is optimal up to a constant factor since every observation must be read, and the sparse algorithm is linear in $T$ and nearly quadratic in $J$. In our experiments, a dense instance with $J=10^5$ and $T=100$, that is, $10^7$ observations, is solved in less than $0.1$ seconds on a laptop, and sparse instances with about $10^3$ dimensions are solved in under a second. Section~\ref{sec:computational_results} describes the experiments, software, and hardware. We have not tried to optimize the running-time for $p\in\{1,\infty\}$, whose purpose is to establish polynomial-time solvability.

\subsection{Outline of the proof of Theorem~\ref{thm:main}}
\label{sec:outline}

The three losses, each in the dense and the sparse setting, give six variants of the problem. All six variants are solved by the same six-step strategy described below. The details of each step depend on the variant. The main idea is that, for a fixed breakpoint $z$, the problem splits into $J$ small fitting problems, one for each dimension. We study how the optimal error of each of these problems varies with $z$, combine these errors over the dimensions, and reduce the search for $z$ to finitely many candidates. 
 
\paragraph{Step 1: reduction to intervals.} The nonlinear terms $(x_i-z)_+$ depend on which observations lie to the right of $z$. Once we restrict $z$ to an interval between consecutive observations, this set is fixed and each residual term, $r_i^j$, becomes affine in $z$. For $i^*\in\{2,\ldots,T-2\}$, let $\textup{OPT}^{i^*}_{p,k}$ denote the optimal value of \eqref{eq:slope-change} with $z\in[x_{i^*},x_{i^*+1}]$. Since these intervals cover $[x_2,x_{T-1}]$, we have that $\textup{OPT}_{p,k}=\min_{i^*\in\{2,\ldots,T-2\}}\textup{OPT}^{i^*}_{p,k}.$ For $z\in[x_{i^*},x_{i^*+1}]$, including the endpoints, the residuals \eqref{eq:residual} become
$$r_i^j=y_i^j-\alpha_j-\beta_jx_i-\delta_jh_i(z),\qquad h_i(z):=\begin{cases}0,& i\leq i^*,\\ x_i-z,& i>i^*.\end{cases}$$
Since $2\leq i^*\leq T-2$, each side of the breakpoint contains at least two observations with distinct time coordinates. We show this implies the three columns $(1)_{i\in[T]}$, $(x_i)_{i\in[T]}$, and $(h_i(z))_{i\in[T]}$ of the fitting problem, corresponding to $\alpha_j$, $\beta_j$, and $\delta_j$, are linearly independent for every $z$, so none of them becomes redundant. We fix $i^*$ for the rest of the outline.
 
\paragraph{Step 2: decoupling for fixed $z$.} For fixed $z$, the dimensions interact only through the choice of which of them may change slope. Each dimension is otherwise fitted separately using the three variables $(\alpha_j,\beta_j,\delta_j)$. For $j\in[J]$ and $z\in[x_{i^*},x_{i^*+1}]$, let $f_j(z)$ be the minimum error when dimension $j$ may change slope at $z$. This error is $\sum_i|r_i^j|$, $\sum_i(r_i^j)^2$, or $\max_i|r_i^j|$ for $p=1,2,\infty$, respectively. Let $U_j$ be the corresponding minimum with $\delta_j=0$. Thus $U_j$ is the error of the best single line, independent of $z$, and $f_j(z)\leq U_j$. Then
$$\textup{OPT}^{i^*}_{p,k}=\min_{z\in[x_{i^*},x_{i^*+1}]}\ \min_{S\subseteq[J],\,|S|\leq k}\ \Big(\bigoplus_{j\in S}f_j(z)\Big)\oplus\Big(\bigoplus_{j\notin S}U_j\Big),$$
where $\oplus$ denotes addition for $p\in\{1,2\}$ and maximum for $p=\infty$, and an empty aggregation equals zero. Thus the problem reduces to understanding and combining the functions $f_j(\cdot)$.
 
\paragraph{Step 3: structure of each $f_j(\cdot)$.} We establish that $f_j(\cdot)$ is continuous. As we cannot evaluate $f_j(z)$ for every $z$, we show that each $f_j(\cdot)$ is piecewise rational on $[x_{i^*},x_{i^*+1}]$, with polynomially many pieces. Each piece has numerator and denominator of degree at most two, and its endpoints are roots of polynomials of degree at most two with rational coefficients.
 
For $p=2$, the least-squares projection gives a single piece, $f_j(z)=U_j-\nu_j(z)^2/D(z)$, where $\nu_j(z)$ is affine in $z$ and $D(z)$ is a positive quadratic in $z$. 
Crucially, $D(\cdot)$ is common to all dimensions.
 
For $p\in\{1,\infty\}$, $f_j(z)$ is the optimal value of a linear program. We show that its feasible region is pointed. Since its objective is bounded below by zero, the linear program has an optimal vertex. 
Since $z$ appears only, and affinely, in the column of $\delta_j$, the Cramer rule shows that every basic solution has entries $(a+bz)/(c+dz)$ with $a,b,c,d\in\Q$. We call such functions linear-fractional. Such a function can be shown to be monotone on every interval that avoids its pole, and two distinct ones cross at most twice. We bound the number of vertices and determine where each is feasible.
\begin{itemize}
\item For $p=\infty$, the linear program minimizes $\eta$ over the four variables $(\alpha_j,\beta_j,\delta_j,\eta)$, subject to the $2T$ constraints $-\eta\leq r_i^j\leq\eta$ for $i\in[T]$. A vertex is determined by four active constraints, so there are $O(T^4)$ candidate vertices. For each candidate, $\eta$ is linear-fractional in $z$, and the candidate is feasible on a union of $O(T)$ intervals. Hence $f_j(\cdot)$ is the minimum of $O(T^5)$ linear-fractional functions, each restricted to one interval.
\item For $p=1$, we prove that an optimal fit can be chosen to pass through three observations, not all on the same side of the breakpoint. Thus, there are $O(T^3)$ choices. Each determines $(\alpha_j,\beta_j,\delta_j)$ uniquely as linear-fractional functions of $z$. Every such fit is feasible, so, unlike for $p=\infty$, no restriction on $z$ is needed. Its loss $\sum_i |r_i^j|$ is piecewise linear-fractional, with $O(T)$ pieces separated by residual sign changes. Thus $f_j(\cdot)$ is the minimum of $O(T^4)$ linear-fractional functions, each restricted to one interval.
\end{itemize}
 
\paragraph{Step 4: aggregation in the dense problem.} Here we minimize $\bigoplus_{j\in[J]}f_j(z)$. Each $f_j$ is piecewise rational on $[x_{i^*},x_{i^*+1}]$, and together their breakpoints partition this interval into cells. On each cell, every $f_j$ is given by a single rational expression whose denominator does not vanish on the cell. We minimize the resulting univariate function on each cell and take the best value over all cells.
\begin{itemize}
\item If $\oplus$ is addition, the objective on each cell is a single rational function. We show that its minimum occurs at a cell endpoint or a root of the numerator of its derivative. This numerator has degree at most two for $p=2$ and at most $2J-2$ for $p=1$.
\item If $\oplus$ is the maximum, the objective on each cell is the upper envelope of $J$ linear-fractional functions. We show each is monotone on the cell, so a minimum occurs at a cell endpoint or a crossing of two functions.
\end{itemize}
 
\paragraph{Step 5: support selection in the sparse problem.} Allowing dimension $j$ to change slope gives the \emph{gain} $u_j(z):=U_j-f_j(z)\geq0$. Our task is to choose at most $k$ dimensions. The optimal choice depends on the aggregation.
 
First, suppose $\oplus$ is addition. For a set $S$, the total error is $\sum_{j\in[J]}U_j-\sum_{j\in S}u_j(z).$ For fixed $z$, an optimal $S$ consists of the $k$ dimensions with the largest gains. The ranking of the gains, however, varies with $z$, and the set of the $k$ largest gains can change only at crossings of two gains, where their relative order may change. Therefore, between consecutive crossings, $S$ is fixed and we solve the dense problem restricted to $S$. For $p=1$, the gains are linear-fractional on each cell of Step~4, so two distinct gains cross at most twice, giving $O(J^2)$ crossings per cell. For $p=2$, the common denominator in $u_j(z)=\nu_j(z)^2/D(z)$ means that the ranking is determined by $|\nu_j(z)|$, which gives the subcells on which the ranking is constant.
 
Next, suppose $\oplus$ is maximum. For $S\subseteq[J]$ with $|S|=k$, the objective is 
\[
\Phi_S(z):=\max \{ \max_{j\in S} f_j(z), \max_{j\notin S}U_j\}.
\]
A dimension outside $S$ keeps its single-line error $U_j$, so the dimensions with the largest $U_j$ should be the ones allowed to change slope. We verify this as follows. Let $S^\star$ be the set of the $k$ dimensions with the largest single-line errors $U_j$, with arbitrary tie-breaking. This set is independent of $z$. We claim that $\Phi_{S^\star}(z)\leq\Phi_S(z)$ for every $z$ and every $S$ of size $k$. If $S\neq S^\star$, choose $j\in S^\star\setminus S$ and $j'\in S\setminus S^\star$. Then $U_j\geq U_{j'}$. Swapping $j$ and $j'$ removes $U_j$ and $f_{j'}(z)$ from the maximum and adds $f_j(z)\leq U_j$ and $U_{j'}\leq U_j$. Thus the maximum cannot increase. Repeated swaps yield $S^\star$, proving the claim. Setting $M:=\max_{j\notin S^\star}U_j$, the problem becomes
$\min_{z\in[x_{i^*},x_{i^*+1}]}\max\Big\{\max_{j\in S^\star}f_j(z),\ M\Big\},$
which is the dense problem on $S^\star$, followed by a maximum with the constant $M$.
 
 
\paragraph{Step 6: candidates and exact comparison.} Steps~4 and~5 reduce the search to finitely many candidates, namely cell endpoints, stationary points, and crossings. These are algebraic numbers of low degree (see Definition~\ref{def:solve}), so we evaluate the objective at each and compare the values exactly. Given an optimal breakpoint $z^\star$, solving the fitting problems of Step~2 at $z^\star$ recovers $(\alpha^\star,\beta^\star,\delta^\star)$.
 
\paragraph{Complexity.} Each $f_j(\cdot)$ has one piece for $p=2$ and $O(T^8)$ pieces for $p\in\{1,\infty\}$. Merging the endpoints of all pieces, and refining by pairwise crossings where needed, gives polynomially many cells. The bounds in Table~\ref{tab:complexity} multiply the number of cells by the cost per cell and by the $T-3$ intervals of Step~1. 

Detailed proofs are presented in Section~\ref{sec:nonsparse} for the dense case and Section~\ref{sec:sparse} for the sparse case.

\section{Polynomial-Time Solvability: Dense Setting}
\label{sec:nonsparse}

This section and the next carry out Steps~3 to~6 of Section~\ref{sec:outline}, with the notation introduced there. Throughout, we fix $i^*\in\{2,\ldots,T-2\}$ and write $\mathcal{I}:=[x_{i^*},x_{i^*+1}]$. By Step~1, the algorithm is run once for each of the $T-3$ values of $i^*$, and the best of the resulting solutions is returned. Hence, all operation counts below are per interval, and they are multiplied by $T$ in Table~\ref{tab:complexity}. For $i\in[T]$ and $z\in\R$, let $a_i(z):=(1,x_i,h_i(z))^\top$, and for $\theta_j=(\alpha_j,\beta_j,\delta_j)^\top$ let $r_i^j(\theta_j,z):=y_i^j-a_i(z)^\top\theta_j$ be the residual of Step~1. For a linear-fractional of the form $(a+bz)/(c+dz)$ with $(c,d)\neq(0,0)$, its derivative is $(bc-ad)/(c+dz)^2$, whose sign is that of the constant $bc-ad$ wherever $c+dz\neq0$. Hence it is strictly increasing, strictly decreasing, or constant on every interval on which $c+dz\neq0$.
\subsection{Preliminaries}

\begin{lemma}\label{lem:basic}
Fix $j\in[J]$ and $p\in\{1,2,\infty\}$.
\begin{enumerate}
\item[(a)] For every $z\in\R$, the vectors $\1$, $x$, and $h(z):=(h_i(z))_{i\in[T]}$ are linearly independent. Equivalently, the vectors $a_i(z)$, $i\in[T]$, span $\R^3$.
\item[(b)] For every $z\in\mathcal{I}$, the minimum defining $f_j(z)$ is attained, and $f_j$ is continuous on $\mathcal{I}$.
\end{enumerate}
\end{lemma}

\begin{proof}
(a) Let $\lambda_1\1+\lambda_2x+\lambda_3h(z)=0$. Since $i^*\geq2$ and $h_1(z)=h_2(z)=0$, we have $\lambda_1+\lambda_2x_i=0$ for $i\in\{1,2\}$, and $x_1\neq x_2$ gives $\lambda_1=\lambda_2=0$. Since $i^*\leq T-2$, we have $h_i(z)=x_i-z$ for $i\in\{T-1,T\}$, and these cannot vanish, so $h(z)\neq0$ and $\lambda_3=0$.

(b) Let $\ell_j(\theta,z)$ denote the loss of dimension $j$, that is, $\sum_i|r_i^j|$, $\sum_i(r_i^j)^2$, or $\max_i|r_i^j|$. We remind the reader that $U_j=\min_{\alpha_j,\beta_j}L_p\big((y_i^j-\alpha_j-\beta_jx_i)_{i\in[T]}\big)$ is the error of the best single line, so $f_j(z)\leq U_j$ for every $z$. Let $\rho:=U_j+1$ for $p\in\{1,\infty\}$ and $\rho:=\sqrt{U_j+1}$ for $p=2$. Let
$$\Theta:=\big\{\theta\in\R^3:\ \ell_j(\theta,z)\leq U_j+1\ \text{for some } z\in\mathcal{I}\big\},$$
which does not depend on $z$. Let $\theta\in\Theta$, and let $z\in\mathcal{I}$ satisfy $\ell_j(\theta,z)\leq U_j+1$. Then $|r_i^j(\theta,z)|\leq\rho$ for every $i$. Using $h_1(z)=h_2(z)=0$, $h_i(z)=x_i-z$ for $i\in\{T-1,T\}$, and
$$r_1^j-r_2^j=y_1^j-y_2^j-\beta(x_1-x_2),\qquad r_T^j-r_{T-1}^j=y_T^j-y_{T-1}^j-(\beta+\delta)(x_T-x_{T-1}),$$
we obtain $|\beta|\leq(|y_1^j-y_2^j|+2\rho)/(x_2-x_1)$ and $|\beta+\delta|\leq(|y_T^j-y_{T-1}^j|+2\rho)/(x_T-x_{T-1})$. Since $\alpha=y_1^j-\beta x_1-r_1^j$, we have $|\alpha|\leq|y_1^j|+|\beta||x_1|+\rho$. Hence $\Theta$ is bounded, and its closure $\bar\Theta$ is compact. 

Fix $z\in\mathcal{I}$. We show that the minimum defining $f_j(z)$ is attained. By definition, $f_j(z)=\inf_{\theta\in\R^3}\ell_j(\theta,z)$, and $f_j(z)\leq U_j<U_j+1$. If $\theta\notin\Theta$, then, by the definition of $\Theta$, $\ell_j(\theta,z')>U_j+1$ for every $z'\in\mathcal{I}$, and in particular $\ell_j(\theta,z)>U_j+1>f_j(z)$. Hence the points outside $\Theta$ do not affect the infimum, and $f_j(z)=\inf_{\theta\in\Theta}\ell_j(\theta,z)$. Moreover, $\Theta\subseteq\bar\Theta\subseteq\R^3$ gives
$$f_j(z)=\inf_{\theta\in\R^3}\ell_j(\theta,z)\leq\inf_{\theta\in\bar\Theta}\ell_j(\theta,z)\leq\inf_{\theta\in\Theta}\ell_j(\theta,z)=f_j(z),$$
so $f_j(z)=\inf_{\theta\in\bar\Theta}\ell_j(\theta,z)$.

Each residual $r_i^j(\theta,z)=y_i^j-\alpha-\beta x_i-\delta h_i(z)$ is continuous in $(\theta,z)$. Since $\ell_j$ is a sum of absolute values ($p=1$), a sum of squares ($p=2$), or a maximum of absolute values ($p=\infty$) of these residuals, $\ell_j$ is continuous on $\R^3\times\mathcal{I}$. In particular, $\ell_j(\cdot,z)$ is continuous, and since $\bar\Theta$ is compact and nonempty, the infimum is attained, that is, $f_j(z)=\min_{\theta\in\bar\Theta}\ell_j(\theta,z)$. As $\ell_j$ is continuous on the compact set $\bar\Theta\times\mathcal{I}$, $f_j(\cdot)$ is continuous on $\mathcal{I}$, by Berge's maximum theorem~\cite{berge1963topological} (Chapter~VI, Section~3). \end{proof}

For $p\in\{1,\infty\}$, $f_j$ is the pointwise minimum of finitely many linear-fractional functions, each available on a part of $\mathcal{I}$. The next lemma turns such a minimum into a piecewise linear-fractional function. It is also used to aggregate the dimensions for $p=\infty$.

\begin{lemma}\label{lem:envelope}
Let $\mathcal{C}=[l,u]$ with $l<u$, where $l$ and $u$ are rational or roots of quadratic polynomials with rational coefficients. For $k\in[m]$, let $\varphi_k=(a_k+b_kz)/(c_k+d_kz)$ have rational coefficients. Let $V_k\subseteq\mathcal{C}$ be the set of points at which $\varphi_k$ takes part in the minimum below. Assume the following.
\begin{enumerate}
\item[(i)] Each $V_k$ is a union of at most $r$ intervals. Each interval may be open, closed, half-open, or a single point.
\item[(ii)] Each endpoint of these intervals is rational, equal to $l$ or $u$, or a root of a quadratic polynomial with rational coefficients.
\item[(iii)] For every $k$, $c_k+d_kz\neq0$ for all $z\in V_k$.
\item[(iv)] $V_1\cup\cdots\cup V_m=\mathcal{C}$, and $\varphi(z):=\min_{k\in[m]:\,z\in V_k}\varphi_k(z)$ is continuous on $\mathcal{C}$.
\end{enumerate}
Then there exist points $l=\zeta_0<\zeta_1<\cdots<\zeta_P=u$ with $P\leq2mr+m^2$, and linear-fractional functions $\psi_1,\ldots,\psi_P$ with rational coefficients, with the following properties. For every $\ell\in[P]$, the denominator of $\psi_\ell$ has no zero in $[\zeta_{\ell-1},\zeta_\ell]$, and $\varphi=\psi_\ell$ on $[\zeta_{\ell-1},\zeta_\ell]$. Each $\zeta_\ell$ is $l$, $u$, an endpoint of an interval of some $V_k$, or a root of a nonzero polynomial of degree at most two with rational coefficients. Given the $\varphi_k$ and the endpoints of the intervals of the $V_k$, the points $\zeta_\ell$ and the functions $\psi_\ell$ are computed with $O((m^2+mr)\log(m + r))$ operations. The same holds with the minimum replaced by the maximum.
\end{lemma}

\begin{proof}
For $k\neq k'$, let $\pi_{kk'}:=(a_k+b_kz)(c_{k'}+d_{k'}z)-(a_{k'}+b_{k'}z)(c_k+d_kz)$, a polynomial of degree at most two. Wherever both functions are defined, we have
$$\varphi_k(z)-\varphi_{k'}(z)=\frac{\pi_{kk'}(z)}{(c_k+d_kz)(c_{k'}+d_{k'}z)}.$$
Let $E$ be the set of points of $(l,u)$ that are endpoints of the intervals forming the sets $V_k$, or roots of a polynomial $\pi_{kk'}\not\equiv0$. Then $|E|\leq2mr+m(m-1)$. Let $\zeta_0<\cdots<\zeta_P$ be the points of $E\cup\{l,u\}$, where points that coincide are listed once, which is detected by exact comparisons as in Remark~\ref{rem:quad-compare}, so that $P\leq|E|+1\leq2mr+m^2$, and consider an open cell $C_\ell:=(\zeta_{\ell-1},\zeta_\ell)$. No endpoint of an interval of $V_k$ lies in $C_\ell$, so either $C_\ell\subseteq V_k$ or $C_\ell\cap V_k=\emptyset$. Let $A_\ell:=\{k:C_\ell\subseteq V_k\}$. Since $V_1\cup\cdots\cup V_m=\mathcal{C}$, the set $A_\ell$ is nonempty.

Let $k,k'\in A_\ell$. Then $C_\ell\subseteq V_k\cap V_{k'}$, so both denominators have no zero in $C_\ell$ by the assumption on $V_k$ and $V_{k'}$. If $\pi_{kk'}\equiv0$, then $\varphi_k=\varphi_{k'}$ on $C_\ell$. Otherwise, $\pi_{kk'}$ has no zero in $C_\ell$, since all its roots in $(l,u)$ belong to $E$. Hence $\varphi_k-\varphi_{k'}$ is continuous and nonzero on $C_\ell$, and by the intermediate value theorem it has a constant nonzero sign on $C_\ell$. Order $A_\ell$ by $k\prec k'$ if $\varphi_k<\varphi_{k'}$ on $C_\ell$, or if $\varphi_k=\varphi_{k'}$ on $C_\ell$ and $k<k'$. This order does not depend on $z\in C_\ell$, so its first element $k_\ell$ satisfies $\varphi=\varphi_{k_\ell}$ on $C_\ell$.

We now extend the formula $\varphi=\varphi_{k_\ell}$ from $C_\ell$ to its closure $[\zeta_{\ell-1},\zeta_\ell]$. If $b_{k_\ell}c_{k_\ell}=a_{k_\ell}d_{k_\ell}$, then $(a_{k_\ell},b_{k_\ell})$ is a multiple of $(c_{k_\ell},d_{k_\ell})$, so $\varphi_{k_\ell}$ is a rational constant on its domain, and we let $\psi_\ell$ be this constant. Otherwise, $\varphi_{k_\ell}$ is nonconstant, and we let $\psi_\ell:=\varphi_{k_\ell}$. The denominator $c_{k_\ell}+d_{k_\ell}z$ of $\psi_\ell$ has no zero in $C_\ell$, since $C_\ell\subseteq V_{k_\ell}$. We show that it has no zero at $\zeta_{\ell-1}$ or $\zeta_\ell$ either. Since $\varphi$ is continuous on the compact interval $[l,u]$, it is bounded. As $\psi_\ell=\varphi$ on $C_\ell$, the function $\psi_\ell$ is bounded on $C_\ell$. A nonconstant linear-fractional function, however, is unbounded near a zero of its denominator. Hence the denominator of $\psi_\ell$ does not vanish at $\zeta_{\ell-1}$ or $\zeta_\ell$. Thus, in both cases, $\psi_\ell$ and $\varphi$ are continuous on $[\zeta_{\ell-1},\zeta_\ell]$ and agree on $C_\ell$, so they agree on $[\zeta_{\ell-1},\zeta_\ell]$.

It remains to compute $k_\ell$ for all $\ell$, and for this we compare neighboring cells. For $\ell\in[P-1]$, call an index $k$ \emph{involved at} $\zeta_\ell$ if $\zeta_\ell$ is an endpoint of an interval of $V_k$ or a root of $\pi_{kk'}\not\equiv0$ for some $k'$, and let $E_{\zeta_\ell}$ be the set of these indices. Indices that are not involved at $\zeta_\ell$ are unaffected when passing from $C_\ell$ to $C_{\ell+1}$. Indeed, let $k\notin E_{\zeta_\ell}$. Since $\zeta_\ell$ is not an endpoint of an interval of $V_k$, either $\zeta_\ell$ is an interior point of an interval of $V_k$, and then $k\in A_\ell\cap A_{\ell+1}$ and $\zeta_\ell\in V_k$, or $k\notin A_\ell\cup A_{\ell+1}$. If, moreover, $k'\notin E_{\zeta_\ell}$ and $k,k'\in A_\ell$, then both denominators are nonzero at $\zeta_\ell\in V_k\cap V_{k'}$, and either $\pi_{kk'}\equiv0$ or $\pi_{kk'}(\zeta_\ell)\neq0$. Hence $\varphi_k-\varphi_{k'}$ is identically zero or has the same sign on $C_\ell$ and $C_{\ell+1}$, so the relative order of $k$ and $k'$ is the same on both cells.

This suggests the following sweep. Sort $E$, which takes $O(|E|\log|E|)$ operations, and visit the cells from left to right, keeping the indices of the current set $A_\ell$ in a balanced binary search tree ordered by $\prec$. When the sweep passes $\zeta_\ell$, delete from the tree the indices of $E_{\zeta_\ell}$ that it contains, and insert those indices of $E_{\zeta_\ell}$ that belong to $A_{\ell+1}$, each at its position on $C_{\ell+1}$. By the observation above, the tree then contains $A_{\ell+1}$ ordered by $\prec$ on $C_{\ell+1}$, and its first element is $k_{\ell+1}$. Each comparison of two indices $k,k'$ during an insertion is decided by the index if $\pi_{kk'}\equiv0$, and otherwise by the signs on $C_{\ell+1}$ of $\pi_{kk'}$ and of the two denominators. Each of these is a nonzero polynomial $g$ of degree at most two with no zero in $C_{\ell+1}$. By the exact Taylor expansion $g(\zeta_\ell+t)=g(\zeta_\ell)+g'(\zeta_\ell)t+\tfrac12g''(\zeta_\ell)t^2$, its sign on $C_{\ell+1}$ is the sign of the first nonzero number among $g(\zeta_\ell)$, $g'(\zeta_\ell)$, and $g''(\zeta_\ell)$. Each index $k$ is involved at no more than the $2r$ endpoints of the intervals of $V_k$ and the $2(m-1)$ roots of the polynomials $\pi_{kk'}$. Hence $\sum_{\ell}|E_{\zeta_\ell}|\leq2mr+2m(m-1)$, and the sweep, including the construction of the initial tree for $C_1$, takes $O((m^2+mr)\log (m + r))$ operations (see \cite{guibas1978dichromatic} for the runtime guarantees of balanced binary search tree operations). For the maximum, the same argument applies with the order of the values reversed.
\end{proof}

\begin{remark}\label{rem:quad-compare}
Two numbers $a+b\sqrt{\Delta}$ and $a'+b'\sqrt{\Delta'}$ with $a,b,a',b',\Delta,\Delta'\in\Q$ and $\Delta,\Delta'\geq0$ can be compared exactly, by sign tests on rational numbers after squaring at most twice, in time polynomial in the encoding lengths of $a,b,a',b',\Delta,\Delta'$.
\end{remark}

\subsection{The $\ell_2$ loss} \label{sec:dense-l2}

Let $A:=[\1\ x]\in\Q^{T\times2}$ and $G:=A^\top A$, which is invertible since the $x_i$ are distinct. Let $M:=I_T-AG^{-1}A^\top$ be the orthogonal projector onto the orthogonal complement of the column space of $A$. Write $h(z)=v-zw$, where $v,w\in\Q^T$ are given by
$$v_i:=\begin{cases}0,& i\leq i^*,\\ x_i,& i>i^*,\end{cases}\qquad w_i:=\begin{cases}0,& i\leq i^*,\\ 1,& i>i^*.\end{cases}$$

Lemma~\ref{lem:l2-closed} establishes the structure of $f_j$, which completes Step~3 for $p=2$.
\begin{lemma}\label{lem:l2-closed}
For every $j\in[J]$ and $z\in\R$,
$$f_j(z)=U_j-\frac{\nu_j(z)^2}{D(z)},\qquad U_j=(y^j)^\top My^j,$$
where
$$\nu_j(z):=(y^j)^\top Mh(z)=s_{j,0}+s_{j,1}z,\qquad D(z):=h(z)^\top Mh(z)=d_0+d_1z+d_2z^2,$$
with $s_{j,0}:=(y^j)^\top Mv$, $s_{j,1}:=-(y^j)^\top Mw$, $d_0:=v^\top Mv$, $d_1:=-2v^\top Mw$, and $d_2:=w^\top Mw$. Moreover, $D(z)>0$ for every $z\in\R$, and the minimum defining $f_j(z)$ is attained at
$$\delta_j=\frac{\nu_j(z)}{D(z)},\qquad(\alpha_j,\beta_j)^\top=G^{-1}A^\top\big(y^j-\delta_jh(z)\big).$$
\end{lemma}

\begin{proof}
For fixed $\delta_j$, the minimum of $\|y^j-\delta_jh(z)-A(\alpha_j,\beta_j)^\top\|_2^2$ over $(\alpha_j,\beta_j)$ is attained at $G^{-1}A^\top(y^j-\delta_jh(z))$, with value $\|M(y^j-\delta_jh(z))\|_2^2$. Since $M$ is symmetric and $M^2=M$, this value equals $U_j-2\delta_j\nu_j(z)+\delta_j^2D(z)$, where $U_j=(y^j)^\top My^j$ is the value for $\delta_j=0$ (that is, the error of the best single line). By Lemma~\ref{lem:basic}(a), $h(z)$ is not in the column space of $A$, so $Mh(z)\neq0$ and $D(z)=\|Mh(z)\|_2^2>0$. Minimizing the quadratic $U_j-2\delta_j\nu_j(z)+\delta_j^2D(z)$ in $\delta_j$ gives the claim.
\end{proof}

\begin{lemma}\label{lem:l2-candidates}
Let $D$ be as in Lemma~\ref{lem:l2-closed}, a polynomial of degree at most two with $D>0$ on $\R$. Let $N$ be a polynomial of degree at most two, let $c\in\R$, and let $F:=c-N/D$. Then $R:=N'D-ND'$ has degree at most two, and $F'=-R/D^2$ on $\R$. Consequently, for every interval $[l,u]$, the minimum of $F$ over $[l,u]$ is attained either at $l$, or at $u$, or, if $R\not\equiv0$, at a root of $R$ in $(l,u)$.
\end{lemma}

\begin{proof}
If $N=n_0+n_1z+n_2z^2$, the coefficient of $z^3$ is $2n_2d_2$ in both $N'D$ and $ND'$, so $\deg R\leq2$. Since $D>0$, $F$ is differentiable on $\R$ with $F'=-R/D^2$. The continuous function $F$ attains its minimum over $[l,u]$. At an interior minimizer, $F'=0$, so $R=0$. If $R\equiv0$, then $F$ is constant, and $l$ is a minimizer.
\end{proof}

\paragraph{Aggregation and candidates (Step 4).} By Step~2, $\textup{OPT}^{i^*}_{2,J}=\min_{z\in\mathcal{I}}F(z)$, where, by Lemma~\ref{lem:l2-closed},
$$F(z):=\sum_{j\in[J]}f_j(z)=C-\frac{N(z)}{D(z)},\qquad C:=\sum_{j\in[J]}U_j,\qquad N(z):=\sum_{j\in[J]}\nu_j(z)^2.$$
The polynomial $N=n_0+n_1z+n_2z^2$, with $n_0=\sum_js_{j,0}^2$, $n_1=2\sum_js_{j,0}s_{j,1}$, and $n_2=\sum_js_{j,1}^2$, has degree at most two irrespective of $J$. By Lemma~\ref{lem:l2-candidates}, an optimal breakpoint on $\mathcal{I}$ lies in the set consisting of $x_{i^*}$, $x_{i^*+1}$, and at most two roots in $(x_{i^*},x_{i^*+1})$ of the quadratic $R=r_0+r_1z+r_2z^2$ with rational coefficients. Hence every candidate lies in $\Q(\sqrt{\Delta})$ with $\Delta:=r_1^2-4r_0r_2\in\Q$, and so does its value $C-N(z)/D(z)$. 
By Remark~\ref{rem:quad-compare}, the values at the at most four candidates are compared exactly in polynomial time, which yields an optimal breakpoint on $\mathcal{I}$.

\begin{remark}
For $J=1$, this procedure is closely related to the classical approach of \cite{hudson1966fitting}, which also examines each interval between consecutive observations. For $J\geq2$, the key point is that $D$ depends only on the time coordinates, so $N$ has degree at most two for every $J$, which is what makes the bound $O(JT)$ possible. Sharing the breakpoint does change the arithmetic, however. For $J=1$, we have $R=\nu_1\cdot(2\nu_1'D-\nu_1D')$, and $2\nu_1'D-\nu_1D'=(2s_{1,1}d_0-s_{1,0}d_1)+(s_{1,1}d_1-2s_{1,0}d_2)z$ is affine. Hence every candidate is rational, and a rational optimal breakpoint exists, whereas for $J=2$ the optimal breakpoint can be irrational (Section~\ref{sec:main}).
\end{remark}

\paragraph{Recovering the solution.} Given an optimal breakpoint $z^\star$, Lemma~\ref{lem:l2-closed} gives $\delta_j^\star=\nu_j(z^\star)/D(z^\star)$ and $(\alpha_j^\star,\beta_j^\star)^\top=G^{-1}A^\top(y^j-\delta_j^\star h(z^\star))$, which lie in $\Q(z^\star)$. This takes $O(JT)$ operations.

\paragraph{Operation count.} Hence the coefficients $d_0,d_1,d_2$ of $D$, the coefficients $s_{j,0},s_{j,1}$ of $\nu_j$, and $U_j$ are determined by $G$, the vectors $A^\top y^j=(\sum_iy_i^j,\sum_ix_iy_i^j)$, the numbers $\sum_i(y_i^j)^2$, and the suffix sums
$$\sum_{i>i^*}x_i,\qquad\sum_{i>i^*}x_i^2,\qquad\sum_{i>i^*}y_i^j,\qquad\sum_{i>i^*}x_iy_i^j,$$
together with $T-i^*$, since $A^\top v=(\sum_{i>i^*}x_i,\sum_{i>i^*}x_i^2)$, $A^\top w=(T-i^*,\sum_{i>i^*}x_i)$, $v^\top v=\sum_{i>i^*}x_i^2$, $v^\top w=\sum_{i>i^*}x_i$, $w^\top w=T-i^*$, $(y^j)^\top v=\sum_{i>i^*}x_iy_i^j$, and $(y^j)^\top w=\sum_{i>i^*}y_i^j$. The quantities that do not depend on $i^*$ are computed once with $O(JT)$ operations. The suffix sums are updated from one value of $i^*$ to the next with $O(J)$ operations, after which $D$, the $s_{j,0},s_{j,1}$, $N$, and $R$ are formed with $O(J)$ operations. Each interval then requires one quadratic root computation and a constant number of comparisons. Over all intervals, and with the recovery step, this gives $O(JT)$ operations.

\subsection{The $\ell_{\infty}$ loss}
\label{sec:dense-linf}

Fix $j$ and we drop the superscript $j$ for simplicity. For $z\in\mathcal{I}$, $f(z)$ is the optimal value of the linear program
\begin{equation}\label{eq:lp-inf}
\min_{\theta\in\R^3,\,\eta\in\R}\ \Big\{\eta:\ \sigma\big(y_i-a_i(z)^\top\theta\big)\leq\eta\ \text{ for all }(i,\sigma)\in[T]\times\{-1,1\}\Big\}.
\end{equation}
Call a set $B\subseteq[T]\times\{-1,1\}$ with $|B|=4$ a \emph{basis}, and let $K_B(z)\in\R^{4\times4}$ be the matrix with rows $(a_i(z)^\top,\sigma)$ for $(i,\sigma)\in B$. The constraint $(i,\sigma)$ holds with equality if and only if $a_i(z)^\top\theta+\sigma\eta=y_i$.

\begin{lemma}\label{lem:linf-bases}
\begin{enumerate}
\item[(a)] For every basis $B$, $\det K_B(z)$ is an affine function of $z$. If $\det K_B\not\equiv0$, let $(\theta_B(z),\eta_B(z))$ be the solution of $a_i(z)^\top\theta+\sigma\eta=y_i$, $(i,\sigma)\in B$. Then $\eta_B$ is linear-fractional with rational coefficients and denominator $\det K_B$, and, for every $(i,\sigma)$, the slack $\eta_B(z)-\sigma r_i(\theta_B(z),z)$ is an affine function with rational coefficients divided by $\det K_B(z)$.
\item[(b)] For each basis $B$ with $\det K_B\not\equiv0$, let 
$$V_B:=\big\{t\in\mathcal{I}:\ \det K_B(t)\neq0\ \text{and all $2T$ slacks of $(\theta_B(t),\eta_B(t))$ are nonnegative}\big\},$$
the set of breakpoints at which $B$ defines a feasible vertex of \eqref{eq:lp-inf}. Each $V_B$ is a union of at most $4T+5$ intervals with rational endpoints. Moreover, for every fixed $z\in\mathcal{I}$,
$$f(z)=\min\big\{\eta_B(z):\ B\text{ a basis with }z\in V_B\big\}.$$\end{enumerate}
\end{lemma}

\begin{proof}
(a) Only the third column of $K_B(z)$ depends on $z$, and its entries $h_i(z)$ are affine in $z$. Expanding along this column shows that $\det K_B(z)$ is affine in $z$. By Cramer's rule, each entry of $(\theta_B(z),\eta_B(z))$ is a determinant divided by $\det K_B(z)$. For $\delta_B$, the third column is replaced by the right-hand side, so the numerator is a rational constant $\kappa$. For $\alpha_B$, $\beta_B$, and $\eta_B$, the numerator still contains the third column, so it is affine. Hence $\eta_B$ is linear-fractional. Moreover, $\delta_B(z)h_i(z)=\kappa h_i(z)/\det K_B(z)$, so every residual $r_i(\theta_B(z),z)=y_i-\alpha_B(z)-\beta_B(z)x_i-\delta_B(z)h_i(z)$, and therefore every slack, is an affine function divided by $\det K_B(z)$.

(b) By Lemma~\ref{lem:basic}(a), the vectors $a_i(z)$ span $\R^3$. Since the rows $(a_i(z)^\top,1)$ and $(-a_i(z)^\top,1)$ of the constraints of \eqref{eq:lp-inf} add up to $(0,0,0,2)$, the constraint matrix has rank four. Hence the feasible region contains no line, that is it is pointed. As \eqref{eq:lp-inf} is feasible and objective bounded below by $0$, it has an optimal vertex \cite{schrijver1998theory}. A vertex is the unique solution of four linearly independent active constraints. Thus it equals $(\theta_B(z),\eta_B(z))$ for a basis $B$ with $\det K_B(z)\neq0$, and then $z\in V_B$. Conversely, for $z\in V_B$, the point $(\theta_B(z),\eta_B(z))$ is feasible, so $\eta_B(z)\geq f(z)$. 

By (a), $\det K_B$ is an affine function with rational coefficients, and it is not identically zero, so it has at most one zero, which is rational. Also by (a), each of the $2T$ slacks has the form $n(t)/\det K_B(t)$, where $n$ is an affine function with rational coefficients. If $n\not\equiv0$, then $n$ has at most one zero, which is rational. If $n\equiv0$, the slack is zero wherever it is defined and never changes sign. Hence the points at which $\det K_B$ or some slack numerator vanishes are at most $2T+1$ rational points. They split $\mathcal{I}$ into at most $2T+2$ open intervals, on each of which all signs are constant, and at most $2T+3$ remaining points, namely these zeros and the two endpoints of $\mathcal{I}$. Hence $V_B$ is a union of at most $(2T+2)+(2T+3)=4T+5$ intervals, some of which may be single points.
\end{proof}

\paragraph{Structure of $f_j$ (Step 3).} We discard the bases with $\det K_B\equiv0$. We apply Lemma~\ref{lem:envelope} on $[l,u]=\mathcal{I}$, whose endpoints are rational, to the functions $\eta_B$ on the sets $V_B$, with $m\leq\binom{2T}{4}=O(T^4)$ and $r\leq4T+5$. Its hypotheses hold for the following reasons. By Lemma~\ref{lem:linf-bases}(a), each $\eta_B$ is linear-fractional with rational coefficients and denominator $\det K_B$, which is nonzero on $V_B$ by the definition of $V_B$. By Lemma~\ref{lem:linf-bases}(b), each $V_B$ is a union of at most $4T+5$ intervals with rational endpoints, the sets $V_B$ cover $\mathcal{I}$, and $f_j(z)=\min_{B:\,z\in V_B}\eta_B(z)$. By Lemma~\ref{lem:basic}(b), $f_j$ is continuous on $\mathcal{I}$. Hence $f_j$ is piecewise linear-fractional with $O(T^8)$ pieces. On each piece, $f_j$ is a linear-fractional function with rational coefficients whose denominator has no zero on the closed piece. The endpoints of the pieces are rational or roots of quadratic polynomials with rational coefficients. Forming the bases and the sets $V_B$ takes $O(T^5\log T)$ operations, and Lemma~\ref{lem:envelope} takes $O(T^8\log T)$ operations.

\paragraph{Aggregation and candidates (Step 4).} By Step~2, $\textup{OPT}^{i^*}_{\infty,J}=\min_{z\in\mathcal{I}}F(z)$ with $F:=\max_{j\in[J]}f_j$. We merge the endpoints of the pieces of all $f_j$, which partitions $\mathcal{I}$ into $O(JT^8)$ closed cells. On a cell $\mathcal{C}$, each $f_j$ is a linear-fractional function $\psi_j$ with rational coefficients whose denominator has no zero on $\mathcal{C}$, and $F$ is continuous. Lemma~\ref{lem:envelope}, applied to the maximum with $m=J$, $r=1$, and $V_j=\mathcal{C}$, splits $\mathcal{C}$ into $O(J^2)$ closed subcells, on each of which $F$ equals one linear-fractional function whose denominator has no zero on the subcell. 

On a subcell $[\alpha,\beta]$, we have $F(z)=\frac{az+b}{cz+d}$ with $cz+d\neq 0$ for all $z\in[\alpha,\beta]$. Since $cz+d$ is continuous and has no zero on $[\alpha,\beta]$, it has constant sign on $[\alpha,\beta]$. Hence, for all $z,w\in[\alpha,\beta]$,
$$ F(z)-F(w)=\frac{(az+b)(cw+d)-(aw+b)(cz+d)}{(cz+d)(cw+d)}=\frac{(ad-bc)(z-w)}{(cz+d)(cw+d)},$$ $$(cz+d)(cw+d)>0.$$
Thus $F$ is nondecreasing on $[\alpha,\beta]$ if $ad-bc\geq 0$ and nonincreasing on $[\alpha,\beta]$ if $ad-bc\leq 0$. In either case,
$$ \min_{z\in[\alpha,\beta]}F(z)=\min\{F(\alpha),F(\beta)\}.$$
Hence an optimal breakpoint is among the endpoints of the subcells, and it remains to select an endpoint with the smallest value of $F$. Each endpoint lies in $\Q(\sqrt{\Delta})$ for some $\Delta\in\Q$, and so does the value of $F$ at it, since $F$ is linear-fractional with rational coefficients on the subcell. By Remark~\ref{rem:quad-compare}, these values are compared exactly, so an optimal breakpoint is found exactly in polynomial time.

\paragraph{Recovering the solution.} Given an optimal breakpoint $z^\star$ and $j\in[J]$, Lemma~\ref{lem:linf-bases}(b) gives an optimal vertex $(\theta_B(z^\star),\eta_B(z^\star))$, where $B$ minimizes $\eta_B(z^\star)$ over the $O(T^4)$ bases with $z^\star\in V_B$. By Cramer's rule, its entries lie in $\Q(z^\star)$. This takes $O(JT^5)$ operations.

\paragraph{Operation count.} For each interval, computing the pieces of all $f_j$ takes $O(JT^8\log T)$ operations, and sorting their endpoints takes $O(JT^8\log(JT))$ operations. Each of the $O(JT^8)$ cells takes $O(J^2\log J)$ operations, and each of the resulting $O(J^3T^8)$ candidates is evaluated with $O(1)$ operations. This gives $O(J^3T^8\log(JT))$ operations per interval, and $O(J^3T^9\log(JT))$ operations in total.

\subsection{The $\ell_1$ loss}
\label{sec:dense-l1}

Fix $j\in[J]$. For convenience, we drop the superscript $j$. For $z\in\mathcal{I}$, $f(z)$ is the optimal value of the linear program
\begin{equation}\label{eq:lp-one}
\min_{\theta\in\R^3,\,e\in\R^T}\ \Big\{\sum_{i\in[T]}e_i:\ \sigma\big(y_i-a_i(z)^\top\theta\big)\leq e_i\ \text{ for all }(i,\sigma)\in[T]\times\{-1,1\}\Big\}.
\end{equation}
For $Q\subseteq[T]$ with $|Q|=3$, let $A_Q(z)\in\R^{3\times3}$ have rows $a_i(z)^\top$, $i\in Q$. If $\det A_Q(z)\neq0$, let $\theta_Q(z):=A_Q(z)^{-1}(y_i)_{i\in Q}$ be the unique fit with zero residual at the observations in $Q$, and let $L_Q(z):=\sum_{i\in[T]}|r_i(\theta_Q(z),z)|$ be its loss.

\begin{lemma}\label{lem:l1-triples}
\begin{enumerate}
\item[(a)] If $\det A_Q(z)\neq0$, then $Q$ contains an index $i\leq i^*$ and an index $i>i^*$.
\item[(b)] For $z\in\mathcal{I}$, $f(z)=\min\{L_Q(z):|Q|=3,\ \det A_Q(z)\neq0\}$.
\item[(c)] If $\det A_Q(z)\not\equiv0$, then $\det A_Q(z) $ is affine in $z$, and $r_i(\theta_Q(z),z)=n_{Q,i}(z)/\det A_Q(z)$ for affine functions $n_{Q,i}$ with rational coefficients. Consequently, $\det A_Q$ has at most one zero, and each of the at most $T-3$ functions $n_{Q,i}\not\equiv0$, $i\notin Q$, has at most one zero. These at most $T-2$ points split $\mathcal{I}$ into at most $T-1$ closed intervals.  On each of them, with the zero of $\det A_Q$ removed, $L_Q=\big(\sum_i\sigma_in_{Q,i}\big)/(\sigma_0\det A_Q)$ for fixed signs $\sigma_0,\sigma_i\in\{-1,1\}$, which is linear-fractional.
\end{enumerate}
\end{lemma}

\begin{proof}
(a) If $Q\subseteq[i^*]$, the third column of $A_Q(z)$ is zero. If $Q\cap[i^*]=\emptyset$, the third column equals the second column minus $z$ times the first. In both cases, $A_Q(z)$ is singular.

(b) We show that $f(z)=\min\big\{L_Q(z):\ Q\subseteq[T],\ |Q|=3,\ \det A_Q(z)\neq0\big\}.$

\emph{Step 1: $f(z)\leq L_Q(z)$ for every such $Q$.} The point $(\theta_Q(z),(|r_i(\theta_Q(z),z)|)_i)$ is feasible for \eqref{eq:lp-one}, and its objective value is $L_Q(z)$. Hence $f(z)\leq L_Q(z)$.

\emph{Step 2: $f(z)=L_Q(z)$ for some such $Q$.} We first show that \eqref{eq:lp-one} has an optimal vertex. Every feasible point satisfies $e_i\geq|r_i(\theta,z)|\geq0$ for all $i$. Let $(\theta,e)+t(\theta',e')$, $t\in\R$, be a line in the feasible region. Since $e+te'\geq0$ for all $t\in\R$, we have $e'=0$. Then $|y_i-a_i(z)^\top\theta-t\,a_i(z)^\top\theta'|\leq e_i$ for all $t\in\R$ and all $i$, so $a_i(z)^\top\theta'=0$ for all $i$. By Lemma~\ref{lem:basic}(a), the vectors $a_i(z)$, $i\in[T]$, span $\R^3$, so $\theta'=0$. Hence the feasible region contains no line. As \eqref{eq:lp-one} is feasible and its objective is bounded below by $0$, it has an optimal vertex $(\bar\theta,\bar e)$ \cite{schrijver1998theory}.

We next count the active constraints at $(\bar\theta,\bar e)$. By optimality, $\bar e_i=|r_i(\bar\theta,z)|$ for every $i$, since otherwise $\bar e_i$ could be decreased. Let $\bar Q:=\{i:r_i(\bar\theta,z)=0\}$. For $i\notin\bar Q$, exactly one of the two constraints of $i$ is active. For $i\in\bar Q$, both are active, and their rows $(a_i(z)^\top,\mathbf{1}_i^\top)$ and $(-a_i(z)^\top,\mathbf{1}_i^\top)$ span the same space as $(a_i(z)^\top,0)$ and $(0,\mathbf{1}_i^\top)$, where $\mathbf{1}_i$ is the $i$th unit vector in $\R^T$. Hence the rank of the active constraints is at most
$$\textup{rank}\{a_i(z):i\in\bar Q\}+|\bar Q|+(T-|\bar Q|)=\textup{rank}\{a_i(z):i\in\bar Q\}+T.$$
At a vertex of a polyhedron in $\R^{T+3}$, this rank equals $T+3$. Therefore the vectors $a_i(z)$, $i\in\bar Q$, span $\R^3$, and some $Q\subseteq\bar Q$ with $|Q|=3$ has $\det A_Q(z)\neq0$. Since $a_i(z)^\top\bar\theta=y_i$ for $i\in Q$, we get $\bar\theta=\theta_Q(z)$, and therefore $L_Q(z)=\sum_i|r_i(\bar\theta,z)|=\sum_i\bar e_i=f(z)$.

Steps 1 and 2 together prove the claim.

(c) As in the proof of Lemma~\ref{lem:linf-bases}(a), $\det A_Q(z)$ is affine, and by Cramer's rule, $\delta_Q$ is a rational constant divided by $\det A_Q(z)$, while $\alpha_Q$ and $\beta_Q$ are affine functions divided by $\det A_Q(z)$. Hence each residual is an affine function $n_{Q,i}$ divided by $\det A_Q(z)$, with $n_{Q,i}\equiv0$ for $i\in Q$. Neither $\det A_Q$ nor any $n_{Q,i}\not\equiv0$ vanishes in the interior of one of the closed intervals. At an endpoint, the $n_{Q,i}$ that vanish contribute zero with either sign. This gives the formula.
\end{proof}

\paragraph{Structure of $f_j$ (Step 3).} By Lemma~\ref{lem:l1-triples}(b) and~(c), $f_j$ is the minimum of the $m\leq\binom{T}{3}(T-1)=O(T^4)$ linear-fractional functions of Lemma~\ref{lem:l1-triples}(c), each available on a single interval. By Lemmas~\ref{lem:basic}(b) and~\ref{lem:envelope} with $r=1$, $f_j$ is piecewise linear-fractional with $O(T^8)$ pieces. On each piece, $f_j$ is a linear-fractional function with rational coefficients whose denominator has no zero on the closed piece. The endpoints of the pieces are rational or roots of quadratic polynomials with rational coefficients. Forming the functions takes $O(T^4\log T)$ operations, and Lemma~\ref{lem:envelope} takes $O(T^8\log T)$ operations.

\paragraph{Aggregation and candidates (Step 4).} By Step~2, $\textup{OPT}^{i^*}_{1,J}=\min_{z\in\mathcal{I}}F(z)$ with $F:=\sum_{j\in[J]}f_j$. Merging the endpoints of the pieces of all $f_j$ partitions $\mathcal{I}$ into $O(JT^8)$ closed cells. On a cell $\mathcal{C}$, $f_j=(a_j+b_jz)/(c_j+d_jz)$ with rational coefficients and $c_j+d_jz\neq0$ on $\mathcal{C}$. Hence $F=P/Q$ on $\mathcal{C}$, where
$$Q(z):=\prod_{j\in[J]}(c_j+d_jz),\qquad P(z):=\sum_{j\in[J]}(a_j+b_jz)\prod_{h\in[J]\setminus\{j\}}(c_h+d_hz),$$
both of degree at most $J$, and $Q\neq0$ on $\mathcal{C}$. Then $F'=H/Q^2$ with $H:=P'Q-PQ'$. If $P=\sum_{t\leq J}p_tz^t$ and $Q=\sum_{t\leq J}q_tz^t$, the coefficient of $z^{2J-1}$ in $H$ is $Jp_Jq_J-Jp_Jq_J=0$, so $\deg H\leq2J-2$. As in the proof of Lemma~\ref{lem:l2-candidates}, the minimum of $F$ over $\mathcal{C}$ is attained at an endpoint of $\mathcal{C}$ or, if $H\not\equiv0$, at a real root of $H$ in the interior of $\mathcal{C}$. We isolate these roots and compare the values of $F$ at all candidates exactly. Every candidate is rational, a root of a quadratic polynomial with rational coefficients, or a root of $H$, which has rational coefficients and degree at most $2J-2$. Hence every candidate and every value of $F$ at a candidate is a real algebraic number of degree at most $2J-2$ whose representation has encoding length polynomial in $T$, $J$, and $\langle\mathcal{D}\rangle$, and such numbers are compared exactly in polynomial time, even when their degrees differ \cite{basu2006algorithms,sagraloff2016computing}.

\paragraph{Recovering the solution.} Given an optimal breakpoint $z^\star$ and $j\in[J]$, Lemma~\ref{lem:l1-triples}(b) gives an optimal fit $\theta_Q(z^\star)$, where $Q$ minimizes $L_Q(z^\star)$ over the triples with $\det A_Q(z^\star)\neq0$. By Cramer's rule, its entries lie in $\Q(z^\star)$. We find such a $Q$ by enumerating the $O(T^3)$ triples and evaluating $L_Q(z^\star)$ for each with $O(T)$ operations, which takes $O(JT^4)$ operations in $\Q(z^\star)$ over all $j$.

\paragraph{Operation count.} For each interval, computing the pieces of all $f_j$ takes $O(JT^8\log T)$ operations, and sorting their endpoints takes $O(JT^8\log(JT))$ operations. On each of the $O(JT^8)$ cells, $P$, $Q$, and $H$ are formed with $O(J^2)$ operations, the real roots of one polynomial of degree at most $2J-2$ are isolated, and at most $2J$ candidates are compared. This gives $O(J^3T^9+JT^9\log(JT))$ arithmetic operations, $O(JT^9)$ root isolations, and $O(J^2T^9)$ algebraic comparisons in total.

\section{Polynomial-Time Solvability: Sparse Setting}
\label{sec:sparse}

We now carry out Step~5 of the proof outline, which handles the sparsity constraint. Let $1\leq k<J$, since the case $k=0$ is solved directly. For $p\in\{1,2\}$, let $u_j:=U_j-f_j\geq0$ be the gains of Step~5. Since the gains are nonnegative, the constraint $|S|\leq k$ may be replaced by $|S|=k$. By Step~2,
$$\textup{OPT}^{i^*}_{p,k}=\min_{z\in\mathcal{I}}\Phi(z),\qquad\Phi(z):=\sum_{j\in[J]}U_j-\max_{S\subseteq[J],\,|S|=k}\ \sum_{j\in S}u_j(z).$$

\begin{lemma}\label{lem:sparse-cells}
Let $p\in\{1,2\}$, let $x_{i^*}=t_0<t_1<\cdots<t_n=x_{i^*+1}$, and for each $q\in[n]$ let $S_q\subseteq[J]$ with $|S_q|=k$ attain the maximum in the definition of $\Phi(z)$ for every $z\in(t_{q-1},t_q)$. Then
$$\textup{OPT}^{i^*}_{p,k}=\min_{q\in[n]}\ \min_{z\in[t_{q-1},t_q]}\Big(\sum_{j\notin S_q}U_j+\sum_{j\in S_q}f_j(z)\Big).$$
\end{lemma}

\begin{proof}
Let $\Phi_q(z):=\sum_{j\notin S_q}U_j+\sum_{j\in S_q}f_j(z)=\sum_{j\in[J]}U_j-\sum_{j\in S_q}u_j(z)$. Then $\Phi_q\geq\Phi$ on $\mathcal{I}$, with equality on $(t_{q-1},t_q)$. Both functions are continuous by Lemma~\ref{lem:basic}(b), so they are equal on $[t_{q-1},t_q]$.
\end{proof}

\subsection{$\ell_2$ loss with sparsity constraint}
\label{sec:sparse-l2}

By Lemma~\ref{lem:l2-closed}, $u_j(z)=U_j-f_j(z)=\nu_j(z)^2/D(z)$, where $\nu_j(z)=(y^j)^\top Mh(z)$ is affine in $z$ and $D(z)=h(z)^\top Mh(z)$ is a quadratic in $z$ that is positive on $\R$ and the same for all $j$. For each $z$, rank the dimensions by $j\prec j'$ if $\nu_j(z)^2>\nu_{j'}(z)^2$, or if $\nu_j(z)^2=\nu_{j'}(z)^2$ and $j<j'$. For every $z$, the first $k$ dimensions in this ranking attain the maximum in the definition of $\Phi(z)$. The ranking can change only at points where $\nu_j(z)^2=\nu_{j'}(z)^2$ for some $j\neq j'$, so we partition $\mathcal{I}$ by these points. Let $E$ be the set of roots in $(x_{i^*},x_{i^*+1})$ of the affine functions $\nu_j(z)-\nu_{j'}(z)$ and $\nu_j(z)+\nu_{j'}(z)$, for the pairs $j\neq j'$ with $\nu_j(z)^2-\nu_{j'}(z)^2\not\equiv0$. These roots are rational, and $|E|\leq J(J-1)$. On each open cell between consecutive points of $E\cup\{x_{i^*},x_{i^*+1}\}$, the function $\nu_j(z)^2-\nu_{j'}(z)^2=(\nu_j(z)-\nu_{j'}(z))(\nu_j(z)+\nu_{j'}(z))$ is either identically zero or has a constant nonzero sign. Hence the ranking, and therefore the set $S$ of its first $k$ dimensions, is constant on the cell. By Lemma~\ref{lem:sparse-cells}, on the closed cell we minimize
$$\sum_{j\in[J]}U_j-\frac{N_S(z)}{D(z)},\qquad N_S(z):=\sum_{j\in S}\nu_j(z)^2,$$
where $N_S(z)$ has degree at most two irrespective of $k$ and $J$. Lemma~\ref{lem:l2-candidates} gives at most four candidates per cell. Each lies in $\Q$ or in a real quadratic field over $\Q$ and, by Remark~\ref{rem:quad-compare}, the candidates are compared exactly in polynomial time. Given an optimal breakpoint $z^\star$ and the set $S$ of its cell, we set $\delta_j^\star=0$ and $(\alpha_j^\star,\beta_j^\star)^\top=G^{-1}A^\top y^j$ for $j\notin S$, and we recover $(\alpha_j^\star,\beta_j^\star,\delta_j^\star)$ for $j\in S$ as in Section~\ref{sec:dense-l2}.

\paragraph{Sweep.} The cells are visited from left to right. On the first cell, the ranking is obtained by sorting, with comparisons at a rational point of the cell. For $t\in E$, let $E_t$ be the set of dimensions $j$ such that $t$ is a root of $\nu_j-\nu_{j'}$ or of $\nu_j+\nu_{j'}$ for some $j'$ with $\nu_j^2-\nu_{j'}^2\not\equiv0$. If $j\notin E_t$, then, for every $j'$, the sign of $\nu_j^2-\nu_{j'}^2$ is the same on the two cells adjacent to $t$, so $j$ has the same rank on both cells. Hence the ranks held by the dimensions in $E_t$ are also the same on both cells. The set $S$ of the cell to the right of $t$ is therefore obtained from the set $S$ of the cell to the left of $t$ by replacing $S\cap E_t$ with the first $|S\cap E_t|$ dimensions of $E_t$ in the ranking on the right cell. This takes $O(|E_t|\log|E_t|)$ operations. The coefficients of $N_S$ are updated with $O(|E_t|)$ operations. Since $\sum_{t\in E}|E_t|\leq2J(J-1)$, the sweep takes $O(J^2\log J)$ operations.

\paragraph{Operation count.} For each interval, the coefficients $s_{j,0},s_{j,1}$ and $d_0,d_1,d_2$ are obtained with $O(J)$ operations from the suffix sums of Section~\ref{sec:dense-l2}. Forming and sorting $E$ takes $O(J^2\log J)$ operations, the sweep takes $O(J^2\log J)$ operations, and each of the $O(J^2)$ cells requires a constant number of operations, quadratic root computations, and comparisons. Over all intervals, and with the recovery step, this gives $O(J^2T\log J)$ operations.

\subsection{$\ell_{\infty}$ loss with sparsity constraint}
\label{sec:sparse-linf}

By Step~5, where $S^\star$ consists of the $k$ dimensions with the largest $U_j$, any minimizer of $g:=\max_{j\in S^\star}f_j$ over $\mathcal{I}$ also minimizes $\max\{g,M\}$, since $t\mapsto\max\{t,M\}$ is nondecreasing. Hence we run the algorithm of Section~\ref{sec:dense-linf} on the dimensions in $S^\star$, which takes $O(k^3T^9\log(kT))$ operations. For $j\notin S^\star$, we use an optimal single line with $\delta_j^\star=0$, obtained as a rational optimal vertex of a linear program with three variables and $2T$ constraints in $O(T)$ operations \cite{megiddo1984linear}. This also gives $U_j$, so computing all $U_j$ and sorting them takes $O(JT+J\log J)$ operations.

\subsection{$\ell_{1}$ loss with sparsity constraint}
\label{sec:sparse-l1}

Each $U_j$ is the optimal value of a linear program, so it can be computed in polynomial time by any polynomial-time linear programming algorithm. To obtain an operation count that does not depend on the encoding length of the data, we instead compute it as follows. The proof of Lemma~\ref{lem:l1-triples}(b), with $a_i(z)$ replaced by $(1,x_i)^\top$, shows that some optimal single line passes through two observations. Hence $U_j$ is the smallest loss among the $\binom{T}{2}$ lines through two observations. Computing all $U_j$, together with optimal lines, takes $O(JT^3)$ operations.

By Section~\ref{sec:dense-l1}, each $f_j$, and hence each gain $u_j=U_j-f_j$, is piecewise linear-fractional with $O(T^8)$ pieces. Merging the endpoints of the pieces of all $f_j$ partitions $\mathcal{I}$ into $O(JT^8)$ closed cells. On a cell $\mathcal{C}$, each $u_j=(a_j+b_jz)/(c_j+d_jz)$ has rational coefficients, and $c_j+d_jz\neq0$ on $\mathcal{C}$. For each pair $j\neq j'$, let $\pi_{jj'}$ be the polynomial of degree at most two defined in the proof of Lemma~\ref{lem:envelope}. We refine $\mathcal{C}$ by the at most $J(J-1)$ roots in its interior of the polynomials $\pi_{jj'}\not\equiv0$. As in the proof of Lemma~\ref{lem:envelope}, on each open subcell the order of the values $u_j$, with ties broken by index, is constant. Hence the set $S$ of the $k$ largest gains is constant on the subcell. We find $S$ with $O(J)$ comparisons by linear-time selection \cite{blum1973time}. To compare $u_j$ and $u_{j'}$ on the subcell, let $\zeta$ be its left endpoint, and note that
$$u_j(z)-u_{j'}(z)=\frac{\pi_{jj'}(z)}{(c_j+d_jz)(c_{j'}+d_{j'}z)}.$$
If $\pi_{jj'}\equiv0$, the tie is broken by index. Otherwise, $\pi_{jj'}$ has no zero on the open subcell, and the denominator has no zero on the closed cell $\mathcal{C}$. Hence the sign of $u_j-u_{j'}$ is constant on the open subcell, and it equals the sign of $(c_j+d_j\zeta)(c_{j'}+d_{j'}\zeta)$ times the sign of $\pi_{jj'}(z)$ for $z>\zeta$ close to $\zeta$. Since $\pi_{jj'}(\zeta+\varepsilon)=\pi_{jj'}(\zeta)+\pi_{jj'}'(\zeta)\varepsilon+\tfrac12\pi_{jj'}''(\zeta)\varepsilon^2$, the latter is the sign of the first nonzero number among $\pi_{jj'}(\zeta)$, $\pi_{jj'}'(\zeta)$, and $\pi_{jj'}''(\zeta)$. As $\zeta$ lies in $\Q$ or in a real quadratic field over $\Q$, these signs are determined exactly with $O(1)$ operations by Remark~\ref{rem:quad-compare}.

 If $\pi_{jj'}\equiv0$, the tie is broken by index. By Lemma~\ref{lem:sparse-cells}, on the closed subcell we minimize $\sum_{j\notin S}U_j+\sum_{j\in S}f_j(z)$. This is the problem of Section~\ref{sec:dense-l1} with $J$ replaced by $k$. Its minimum is attained at an endpoint of the subcell or at a root of a polynomial $H_S$ of degree at most $2k-2$, formed with $O(k^2)$ operations. For $k=1$, $H_S$ is constant, and only the endpoints are candidates. Given an optimal breakpoint $z^\star$ and the set $S$ of its subcell, the fits of the dimensions in $S$ are recovered as in Section~\ref{sec:dense-l1}, and an optimal line is used for every $j\notin S$.

\paragraph{Operation count.} For each interval, computing the pieces of all $f_j$ and sorting their endpoints takes $O(JT^8\log(JT))$ operations. Each of the $O(JT^8)$ cells requires $O(J^2\log J)$ operations to form and sort the crossings, and each of its $O(J^2)$ subcells requires $O(J+k^2)$ operations, one root isolation of a polynomial of degree at most $2k-2$, and $O(k)$ algebraic comparisons. Together with the computation of the $U_j$, this gives $O((J+k^2)J^3T^9+JT^9\log(JT))$ arithmetic operations, $O(J^3T^9)$ root isolations, and $O(kJ^3T^9)$ algebraic comparisons in total.

\section{Computational Results}
\label{sec:computational_results}


In this section, we evaluate the practical performance of our algorithms for the least-squares loss ($p=2$). We solve the dense and sparse problems on simulated time series and on three public datasets, and we compare our algorithms against Gurobi applied to a mixed-integer nonlinear programming (MINLP) formulation.
\subsection{Scope of the experiments}
\label{sec:comp-scope}

\paragraph{Why only the $\ell_2$ loss.} We restrict the experiments to $p=2$. Least squares is the standard loss for fitting a continuous piecewise-linear trend, and it corresponds to maximum likelihood estimation under Gaussian noise \cite{hudson1966fitting,kim2011estimating,fearnhead2019detecting}. It is also the case in which our algorithms are practical. By Table~\ref{tab:complexity}, the dense algorithm is linear in the size of the input, and the sparse algorithm performs $O(J^2T\log J)$ operations. Although the algorithms for $p\in\{1,\infty\}$ run in polynomial time, their high-degree dependence on $T$ makes them impractical for moderate-sized instances. This is not merely an artifact of the worst-case analysis: for $p=\infty$, the algorithm explicitly enumerates all $\binom{2T}{4}$ bases of \eqref{eq:lp-inf} for each dimension and interval, resulting in more than $6.4\times 10^7$ bases per dimension when $T=100$.

\paragraph{Why floating-point arithmetic.} Theorem~\ref{thm:main} guarantees exact computation, but our implementation uses double-precision arithmetic, so its solutions are approximately optimal and not certified in the sense of Definition~\ref{def:solve}. Exact arithmetic is possible, but for large $J$ its cost would be dominated by large-integer operations, since the coefficients of $N$ aggregate all $J$ dimensions.

\subsection{Datasets}
\label{sec:comp-datasets}

Table~\ref{tab:dataset_design} summarizes the instances. Each dataset
gives one dense instance and two sparse instances with
$k=\lfloor0.20J\rfloor$ and $k=\lfloor0.10J\rfloor$.
We consider sparsity levels of $20\%$ and $10\%$ to evaluate the
algorithms under different degrees of sparsity, where changes are
restricted to a subset of dimensions, motivated by similar sparsity
levels considered in the computational experiments of
\cite{hahn2020bayesproject,anastasiou2023generalized}. The subsets of the three public datasets used in our experiments, together with visualizations, are available in the repositories listed below.


\paragraph{Simulated data.}
We generate correlated stock-price trajectories with
market, sector, and stock-specific fluctuations, as
described in Appendix~\ref{app:synthetic_data}.
Each dimension independently undergoes a change in
its expected daily log-return with a fixed probability,
at a randomly selected time near the middle of the
observation horizon. We use $T=100$ and $x_i=i-1$.

\paragraph{ERA5 temperature data.} ERA5 is a global atmospheric reanalysis produced by the European Centre for Medium-Range Weather Forecasts, which combines historical observations with numerical weather models \cite{hersbach2020era5}. Motivated by studies of changing warming rates \cite{cahill2015change,rahmstorf2017global,foster2026global}, we use the $2$\,m air temperature on a $201\times201$ latitude-longitude grid with $0.25^\circ$ spacing, covering $10^\circ$N to $60^\circ$N and $125^\circ$W to $75^\circ$W, including land and ocean. Each grid point is one dimension, so $J=40{,}401$. We use annual values for 2000 to 2025, so $T=26$ and $x_i=1999+i$, where each value is the average over the year of the daily temperatures at 17:00 UTC. The data are available at \url{https://github.com/ApoorvaGOR/ERA5-Regional-Temperature}.

\paragraph{HoloAssist hand-tracking data.} The HoloAssist dataset \cite{wang2023holoassist} contains hand movements tracked by the Microsoft HoloLens mixed-reality headset. We use an $8.013$-second sequence sampled at about $30$ frames per second, so $T=240$, and the time coordinates are the timestamps in seconds. Each hand has $26$ tracked joints with $21$ measurements per joint, namely $16$ entries of a pose transformation matrix, $3$ spatial coordinates, and $2$ tracking-status indicators. This gives $546$ dimensions per hand and $J=1{,}092$ for both hands. The data are available at \url{https://github.com/ApoorvaGOR/HoloAssist-Hand-Sample}.

\paragraph{NASA C-MAPSS turbofan data.} The Commercial Modular Aero-Propulsion System Simulation (C-MAPSS) dataset \cite{saxena2008turbofan} contains multivariate sensor measurements from simulated aircraft engines undergoing progressive degradation, a setting that motivates multidimensional slope change detection \cite{cao2018multi}. We use training engine~1 of the FD001 subset, with all $21$ sensor channels over its complete trajectory of $192$ cycles, so $J=21$ and $T=192$. Operational settings and remaining useful life (RUL) labels are not fitted. With $R_0=191$ the initial RUL and $R_i$ the RUL at observation $i$, we set $x_i=R_0+1-R_i$, which gives the operating cycles $1,\ldots,192$. The data are available at \url{https://github.com/ApoorvaGOR/C-MAPSS-Turbofan-Data}.

\begin{table}[!ht]
\centering
\caption{Instances considered in the computational experiments, including the dense setting and the sparse setting with sparsity levels of $20\%$ and $10\%$.}
\label{tab:dataset_design}

\scriptsize
\renewcommand{\arraystretch}{1.2}
\setlength{\tabcolsep}{3pt}

\begin{tabular}{l|l|l|l|c|l}
\hline

\textbf{Dataset}
& \textbf{Dense: $J$}
& \textbf{Sparse: $(J;k)$, $20\%$}
& \textbf{Sparse: $(J;k)$, $10\%$}
& \textbf{$T$}
& \textbf{Time coordinates} \\
\hline

Simulated
& \begin{tabular}{@{}l@{}}
$10,20,10^2,10^3,$\\
$10^4,10^5$
\end{tabular}
& \begin{tabular}{@{}l@{}}
$(10;2),(20;4),(10^2;20)$\\
$(10^3;200),(10^4;2000)$\\
$(10^5;20000)$
\end{tabular}
& \begin{tabular}{@{}l@{}}
$(10;1),(20;2),(10^2;10)$\\
$(10^3;100),(10^4;1000)$\\
$(10^5;10000)$
\end{tabular}
& 100
& $0,1,\ldots,99$ \\
\hline

ERA5
& $40{,}401$
& $(40{,}401;8{,}080)$
& $(40{,}401;4{,}040)$
& 26
& Years 2000 to 2025 \\
\hline

HoloAssist
& $1{,}092$
& $(1{,}092;218)$
& $(1{,}092;109)$
& 240
& Timestamps (s) \\
\hline

NASA C-MAPSS
& $21$
& $(21;4)$
& $(21;2)$
& 192
& Cycles $1,\ldots,192$ \\
\hline

\end{tabular}
\end{table}



\subsection{Implementation and benchmark}
\label{sec:comp-implementation}


\paragraph{Implementation and computing environment.} The algorithms are implemented in Python, and their computationally intensive routines are compiled with the just-in-time (JIT) compiler Numba \cite{lam2015numba}. The experiments were run on a MacBook Pro (November 2023) with an Apple M3 chip and 16~GB of unified memory, running macOS Golden Gate~27.0.1, and the MINLP was solved with Gurobi Optimizer 13.0.2 \cite{gurobi}. Both our algorithms
and the Gurobi MINLP benchmark were run using a single
CPU thread. For our algorithms, we restricted Numba and
the underlying numerical libraries to one thread, while
Gurobi was configured to use one thread. The source code and computational results are publicly available at the following repository:

\noindent\url{https://github.com/ApoorvaGOR/High-Dimensional-Piecewise-Linear-Fitting}.

\paragraph{Runtime measurement.}
Numba compiles each numerical routine the first time it is called.
Before the experiments, we therefore run each algorithm once on a small instance,
with $J=2$ and $T=8$ for the dense algorithm and with $J=4$, $T=8$, and $k=2$
for the sparse algorithm.
The experiments for the two sparsity levels ($20\%$ and $10\%$) are conducted
in separate executions, each with its own warmup.
The warmup is performed once per execution, and the compiled code is then
reused for all instances within that execution.
The warmup takes $3.613$\,s for the dense algorithm and $5.090$\,s and
$5.120$\,s for the sparse algorithm in the $20\%$ and $10\%$ sparsity
experiments, respectively.
The reported runtimes of our algorithms exclude the warmup,
data acquisition or generation, time-coordinate conversion, and output storage.
All runtimes are from single runs.
Runtimes are measured using Python's \texttt{time.perf\_counter()},
a high-resolution monotonic wall-clock timer, which uses macOS's
\texttt{mach\_absolute\_time()} with a reported resolution of
approximately $41.67$\,ns on our machine.


\paragraph{MINLP benchmark.} As a benchmark, we use the MINLP formulation of \cite{narula2026changepointdetectionpiecewiselinear} for a single shared breakpoint, stated with its big-$M$ constants in Appendix~\ref{app:minlp_bigM}. It has $T$ binary variables that assign the observations to the two segments, $J$ further binary variables in the sparse problem, and $J$ bilinear continuity constraints. Each instance is solved by Gurobi with a time limit of $10,000$\,s. The reported MINLP runtime is the solver time, and it excludes the construction of the parameter bounds and of the model. Gurobi may exceed the time limit slightly while terminating. We report the best lower bound (LB), the objective value of the best solution found (UB), and the relative optimality gap. Our algorithms are run without a time limit.


\subsection{Results}
\label{sec:comp-results}

Table~\ref{tab:all_results} reports the results for the dense and the sparse problem. The column $OPT_2$ gives the optimal value, that is, the sum of squared errors, computed by our algorithms.

\begin{table}[!ht]
\centering
\caption{Computational comparison of our algorithms and the MINLP benchmark
for (a) the dense problem of Section~\ref{sec:dense-l2},
(b) the sparse problem of Section~\ref{sec:sparse} with
$k=\lfloor 0.20J\rfloor$, and
(c) the sparse problem with $k=\lfloor 0.10J\rfloor$.}
\label{tab:all_results}

\small
\setlength{\tabcolsep}{4pt}

\resizebox{\textwidth}{!}{%
\begin{tabular}{l r r r | r r | r r r r}
\hline
& & & &
\multicolumn{2}{c|}{Our algorithm} &
\multicolumn{4}{c}{MINLP benchmark} \\

Dataset & $J$ & $T$ & $k$ &
$OPT_2$ & Runtime (s) &
LB & UB & Gap (\%) & Runtime (s) \\
\hline

\multicolumn{10}{l}{\textbf{(a) Dense problem}} \\
\hline

Simulated
& 10 & 100 & --
& \textbf{25,830.48}
& $\mathbf{2.600\times 10^{-4}}$
& 25,827.92
& \textbf{25,830.48}
& 0.01
& $1.726\times 10^{1}$ \\

& 20 & 100 & --
& \textbf{78,409.82}
& $\mathbf{2.020\times 10^{-4}}$
& 78,403.28
& \textbf{78,409.82}
& 0.01
& $5.374\times 10^{1}$ \\

& 100 & 100 & --
& \textbf{417,949.81}
& $\mathbf{2.090\times 10^{-4}}$
& 417,908.01
& \textbf{417,949.81}
& 0.01
& $7.097\times 10^{2}$ \\

& 1,000 & 100 & --
& \textbf{3,691,434.08}
& $\mathbf{4.790\times 10^{-4}}$
& 3,683,252.43
& \textbf{3,691,434.08}
& 0.22
& $1.000\times 10^{4}$ $^{*}$ \\

& 10,000 & 100 & --
& \textbf{36,169,676.73}
& $\mathbf{4.769\times 10^{-3}}$
& 17,066,148.10
& 40,817,159.15
& 58.19
& $1.000\times 10^{4}$ $^{*}$ \\

& 100,000 & 100 & --
& \textbf{358,712,206.42}
& $\mathbf{7.961\times 10^{-2}}$
& 0.04
& 607,422,150.29
& 100.00
& $1.002\times 10^{4}$ $^{*}$ \\

\hline

ERA5
& 40,401 & 26 & --
& \textbf{501,263.44}
& $\mathbf{4.681\times 10^{-3}}$
& $8.55\times 10^{-6}$
& 528,198.01
& 100.00
& $1.000\times 10^{4}$ $^{*}$ \\

HoloAssist
& 1,092 & 240 & --
& \textbf{6,176.35}
& $\mathbf{3.274\times 10^{-3}}$
& 5,878.88
& 6,176.43
& 4.82
& $1.000\times 10^{4}$ $^{*}$ \\

NASA C-MAPSS
& 21 & 192 & --
& \textbf{11,686.49}
& $\mathbf{2.960\times 10^{-4}}$
& 11,685.34
& \textbf{11,686.49}
& 0.01
& $7.520\times 10^{0}$ \\

\hline
\multicolumn{10}{l}{\textbf{(b) Sparse problem ($k=\lfloor0.20J\rfloor$)}} \\
\hline

Simulated
& 10 & 100 & 2
& \textbf{31,428.50}
& $\mathbf{1.950\times 10^{-4}}$
& 31,425.62
& \textbf{31,428.50}
& 0.01
& $4.090\times 10^{0}$ \\

& 20 & 100 & 4
& \textbf{91,112.37}
& $\mathbf{4.270\times 10^{-4}}$
& 91,103.76
& \textbf{91,112.37}
& 0.01
& $9.139\times 10^{1}$ \\

& 100 & 100 & 20
& \textbf{500,700.13}
& $\mathbf{6.138\times 10^{-3}}$
& 369,423.81
& 503,853.64
& 26.68
& $1.000\times 10^{4}$ $^{*}$ \\

& 1,000 & 100 & 200
& \textbf{4,381,778.86}
& $\mathbf{5.797\times 10^{-1}}$
& 2,077,935.04
& 4,399,461.26
& 52.77
& $1.000\times 10^{4}$ $^{*}$ \\

& 10,000 & 100 & 2,000
& \textbf{43,815,769.93}
& $\mathbf{6.768\times 10^{1}}$
& $1.09\times10^{-5}$
& 53,566,082.59
& 100.00
& $1.000\times 10^{4}$ $^{*}$ \\

& 100,000 & 100 & 20,000
& --
& --$^{\dagger}$
& 0.04
& 617,148,417.09
& 100.00
& $1.002\times 10^{4}$ $^{*}$ \\

\hline

ERA5
& 40,401 & 26 & 8,080
& \textbf{513,464.42}
& $\mathbf{8.881\times 10^{3}}$
& $8.55\times10^{-6}$
& 537,144.62
& 100.00
& $1.000\times 10^{4}$ $^{*}$ \\

HoloAssist
& 1,092 & 240 & 218
& \textbf{6,485.05}
& $\mathbf{7.153\times 10^{-1}}$
& $1.60\times10^{-8}$
& 7,921.85
& 100.00
& $1.000\times 10^{4}$ $^{*}$ \\

NASA C-MAPSS
& 21 & 192 & 4
& \textbf{11,766.44}
& $\mathbf{6.300\times 10^{-4}}$
& 11,765.34
& \textbf{11,766.44}
& 0.01
& $1.216\times 10^{1}$ \\

\hline
\multicolumn{10}{l}{\textbf{(c) Sparse problem ($k=\lfloor0.10J\rfloor$)}} \\
\hline

Simulated
& 10 & 100 & 1
& \textbf{37,049.63}
& $\mathbf{1.710\times 10^{-4}}$
& 37,049.63
& 37,051.92
& 0.01
& $1.638\times 10^{0}$ \\

& 20 & 100 & 2
& \textbf{105,367.14}
& $\mathbf{3.960\times 10^{-4}}$
& 105,362.30
& \textbf{105,367.14}
& 0.00
& $1.064\times 10^{1}$ \\

& 100 & 100 & 10
& \textbf{558,320.79}
& $\mathbf{6.381\times 10^{-3}}$
& 386,192.66
& 559,523.91
& 30.98
& $1.000\times 10^{4}$ $^{*}$ \\

& 1,000 & 100 & 100
& \textbf{4,888,999.40}
& $\mathbf{6.119\times 10^{-1}}$
& 1,063,321.38
& 5,225,165.60
& 79.65
& $1.000\times 10^{4}$ $^{*}$ \\

& 10,000 & 100 & 1,000
& \textbf{49,101,280.35}
& $\mathbf{6.590\times 10^{1}}$
& $1.09\times 10^{-5}$
& 62,069,667.92
& 100.00
& $1.000\times 10^{4}$ $^{*}$ \\

& 100,000 & 100 & 10,000
& $-$
& $-$ $^{\dagger}$
& 0.04
& 617,148,417.09
& 100.00
& $1.002\times 10^{4}$ $^{*}$ \\

\hline

ERA5
& 40,401 & 26 & 4,040
& \textbf{521,415.53}
& $\mathbf{7.704\times 10^{3}}$
& $8.55\times 10^{-6}$
& 537,144.62
& 100.00
& $1.000\times 10^{4}$ $^{*}$ \\

HoloAssist
& 1,092 & 240 & 109
& \textbf{6,807.89}
& $\mathbf{7.194\times 10^{-1}}$
& $1.60\times 10^{-8}$
& 8,179.81
& 100.00
& $1.000\times 10^{4}$ $^{*}$ \\

NASA C-MAPSS
& 21 & 192 & 2
& \textbf{12,756.23}
& $\mathbf{1.181\times 10^{-3}}$
& 12,755.38
& \textbf{12,756.23}
& 0.01
& $1.048\times 10^{1}$ \\

\hline
\end{tabular}%
}

\vspace{0.3em}

\begin{minipage}{\textwidth}
\footnotesize
\textit{Notes:}
The one-time JIT compilation times of approximately
$3.613$\,s for the dense algorithm,
$5.090$\,s for the sparse algorithm with $20\%$ sparsity,
and $5.120$\,s for the sparse algorithm with $10\%$ sparsity
are excluded from the reported runtimes.
LB, UB, and the gap refer to the bounded MINLP, with the gap
representing Gurobi's relative gap.
In each row, bold runtimes indicate the shorter runtime,
and bold objective values indicate the lowest displayed
objective value.
$^{*}$The MINLP reached the prescribed time limit of
$10,000$\,s.
$-^{\dagger}$Our algorithm did not terminate within the
prescribed time limit of $10,000$\,s; consequently, its
runtime and objective value are not reported.
\end{minipage}

\end{table}

\paragraph{Dense problem.} The dense algorithm solves every instance in less than $0.1$\,s, including the simulated instance with $J=10^5$ and $10^7$ observations, which takes less than $0.080$\,s. This speed matters in online detection \cite{cao2018multi}, where the problem is solved again after each new observation, for example to evaluate the generalized likelihood ratio statistic \cite{siegmund1995using}, which under Gaussian noise is proportional to $\sum_{j\in[J]}U_j-\textup{OPT}_{2,J}$. On HoloAssist, one solve takes $3.27$\,ms, well below the interval of about $33$\,ms between consecutive frames at $30$ frames per second.


\paragraph{Sparse problem.}
For both sparsity levels ($20\%$ and $10\%$), the sparse algorithm solves all instances with $J\leq1092$ in less than $1$\,s. For the simulated instance with $J=10^4$, the runtimes are $67.68$\,s and $65.90$\,s for $20\%$ and $10\%$ sparsity, respectively. The ERA5 instance, with $J=40{,}401$, takes approximately $2.5$ hours and $2.1$ hours, respectively. For the simulated instance with $J=10^5$, the algorithm does not finish within $10{,}000$\,s at either sparsity level. These results reflect the quadratic dependence on $J$ in the bound $O(J^2T\log J)$, since each interval can contain up to $J(J-1)$ crossings of two gains, approximately $1.6\times10^9$ for ERA5. Hence, the sparse algorithm remains practical for instances with up to about $10^4$ dimensions at both sparsity levels.

\paragraph{Comparison with the MINLP.}
Gurobi solves the MINLP to its default relative gap of $0.01\%$
only on the NASA C-MAPSS instances and on the simulated instances
with $J\leq 100$ in the dense problem and $J\leq 20$ in the
sparse problem, at both $20\%$ and $10\%$ sparsity levels.
On these instances, our algorithms obtain objective values no worse
than Gurobi's, and are faster by factors ranging from approximately
$8.9\times 10^3$ to $3.4\times 10^6$.
On all other instances, Gurobi reaches the time limit of
$10{,}000$\,s with gaps between $0.22\%$ and $100\%$.
Its best solution matches ours only on the dense simulated instance
with $J=1000$ and is otherwise worse, for example by $5.4\%$
on the dense ERA5 instance and by $69\%$ on the dense simulated
instance with $J=10^5$.
For the sparse simulated instance with $J=10^5$, neither method
finishes within $10{,}000$\,s at either sparsity level.

\section{Conclusion}
\label{sec:conclusion}

We studied the problem of fitting continuous piecewise-linear functions with a single breakpoint shared by all dimensions of a multidimensional time series, where the breakpoint varies continuously. We showed that this problem is solvable in polynomial time under the $\ell_1$, $\ell_2$, and $\ell_\infty$ losses, both in the dense setting and in the sparse setting, in which at most $k$ dimensions may change slope. Since optimal breakpoints can be irrational, we specified what it means to solve the problem exactly. Our algorithms return an optimal breakpoint as a real algebraic number, of degree at most two for $p\in\{2,\infty\}$ and of degree bounded in terms of $J$ or $k$ for $p=1$, and they return the fitting coefficients as rational functions of it. All six variants are handled by the same strategy, which involves reduction to intervals between consecutive observations, decoupling of the dimensions for a fixed breakpoint, a description of each dimension-wise optimal error as a function of the breakpoint, and aggregation of these functions over the dimensions.

For the least-squares loss, the algorithms are also practical. On simulated data and on climate, hand-tracking, and turbofan data, the dense algorithm solved every instance, with up to $10^7$ observations, in less than $0.1$\,s, and the sparse algorithm solved instances with about a thousand dimensions in less than a second. In particular, the dense algorithm is fast enough to solve the problem again each time a new observation arrives, which makes exact least-squares fitting usable for online slope change detection. A state-of-the-art MINLP solver, applied to an exact mixed-integer formulation, solved only the smallest instances within $10,000$\,s and otherwise ended with large optimality gaps.

Several questions remain open. The sparse algorithm took approximately $2.5$ hours and $2.1$ hours on the ERA5 instance with $40{,}401$ dimensions at $20\%$ and $10\%$ sparsity, respectively, primarily due to the large number of pairwise crossings of the gains that must be processed. An algorithm that tracks only the crossings that change the set of the $k$ largest gains could improve this substantially. For $p\in\{1,\infty\}$, our operation counts are polynomial but of high degree in $T$, and it would be valuable to find algorithms whose dependence on $T$ is close to that of the least-squares case. Finally, the complexity of fitting several shared breakpoints remains open. This problem is equivalent to globally training a ReLU network with time as its input, several hidden units, and $J$ outputs that share those units. It is natural to ask whether the problem stays polynomial-time solvable when the number of breakpoints is fixed, and whether it becomes NP-hard when that number is part of the input.

\newpage

\section*{Acknowledgements}
The authors acknowledge the use of large language models to assist in writing the code for the computational experiments in Section~\ref{sec:computational_results}. All code was reviewed and tested by the authors, who take full responsibility for its correctness and for the reported results.

\bibliographystyle{plain}
\bibliography{refs}

@article{toriello2012fitting,
  author  = {Toriello, Alejandro and Vielma, Juan Pablo},
  title   = {Fitting piecewise linear continuous functions},
  journal = {European Journal of Operational Research},
  volume  = {219},
  number  = {1},
  pages   = {86--95},
  year    = {2012}
}

@article{blum1973time,
  author  = {Blum, Manuel and Floyd, Robert W. and Pratt, Vaughan and Rivest, Ronald L. and Tarjan, Robert E.},
  title   = {Time bounds for selection},
  journal = {Journal of Computer and System Sciences},
  volume  = {7},
  number  = {4},
  pages   = {448--461},
  year    = {1973},
  doi     = {10.1016/S0022-0000(73)80033-9}
}

@book{berge1963topological,
  author    = {Berge, Claude},
  title     = {Topological Spaces},
  publisher = {Oliver and Boyd},
  address   = {Edinburgh and London},
  year      = {1963}
}

@article{bienstock2023complexity,
  author  = {Bienstock, Daniel and Del Pia, Alberto and Hildebrand, Robert},
  title   = {Complexity, Exactness, and Rationality in Polynomial Optimization},
  journal = {Mathematical Programming},
  volume  = {197},
  number  = {2},
  pages   = {661--692},
  year    = {2023}
}

@article{delpia2017mixed,
  author  = {Del Pia, Alberto and Dey, Santanu S. and Molinaro, Marco},
  title   = {Mixed-Integer Quadratic Programming is in {NP}},
  journal = {Mathematical Programming},
  volume  = {162},
  number  = {1--2},
  pages   = {225--240},
  year    = {2017}
}

@book{basu2006algorithms,
  author    = {Basu, Saugata and Pollack, Richard and Roy, Marie-Fran{\c{c}}oise},
  title     = {Algorithms in Real Algebraic Geometry},
  edition   = {Second},
  series    = {Algorithms and Computation in Mathematics},
  volume    = {10},
  publisher = {Springer},
  address   = {Berlin},
  year      = {2006}
}

@article{megiddo1984linear,
  author  = {Megiddo, Nimrod},
  title   = {Linear Programming in Linear Time When the Dimension Is Fixed},
  journal = {Journal of the ACM},
  volume  = {31},
  number  = {1},
  pages   = {114--127},
  year    = {1984}
}

@article{jackson2005algorithm,
  author  = {Jackson, Brad and Scargle, Jeffrey D. and Barnes, David and Arabhi, Sundararajan and Alt, Alina and Gioumousis, Peter and Gwin, Elyus and Sangtrakulcharoen, Paungkaew and Tan, Linda and Tsai, Tun Tao},
  title   = {An Algorithm for Optimal Partitioning of Data on an Interval},
  journal = {IEEE Signal Processing Letters},
  volume  = {12},
  number  = {2},
  pages   = {105--108},
  year    = {2005}
}

@article{killick2012optimal,
  author  = {Killick, Rebecca and Fearnhead, Paul and Eckley, Idris A.},
  title   = {Optimal Detection of Changepoints with a Linear Computational Cost},
  journal = {Journal of the American Statistical Association},
  volume  = {107},
  number  = {500},
  pages   = {1590--1598},
  year    = {2012}
}

@article{kim2009trend,
  author  = {Kim, Seung-Jean and Koh, Kwangmoo and Boyd, Stephen and Gorinevsky, Dimitry},
  title   = {$\ell_1$ Trend Filtering},
  journal = {SIAM Review},
  volume  = {51},
  number  = {2},
  pages   = {339--360},
  year    = {2009}
}

@article{tibshirani2014adaptive,
  author  = {Tibshirani, Ryan J.},
  title   = {Adaptive Piecewise Polynomial Estimation via Trend Filtering},
  journal = {The Annals of Statistics},
  volume  = {42},
  number  = {1},
  pages   = {285--323},
  year    = {2014}
}

@article{hudson1966fitting,
  author  = {Hudson, Derek J.},
  title   = {Fitting Segmented Curves Whose Join Points Have to Be Estimated},
  journal = {Journal of the American Statistical Association},
  volume  = {61},
  number  = {316},
  pages   = {1097--1129},
  year    = {1966}
}

@article{muggeo2003estimating,
  author  = {Muggeo, Vito M. R.},
  title   = {Estimating Regression Models with Unknown Break-Points},
  journal = {Statistics in Medicine},
  volume  = {22},
  number  = {19},
  pages   = {3055--3071},
  year    = {2003}
}

@misc{bleakley2011group,
  author        = {Bleakley, Kevin and Vert, Jean-Philippe},
  title         = {The Group Fused Lasso for Multiple Change-Point Detection},
  year          = {2011},
  eprint        = {1106.4199},
  archivePrefix = {arXiv},
  note = {\url{https://arxiv.org/abs/1106.4199}}
}

@inproceedings{arora2018understanding,
  author    = {Arora, Raman and Basu, Amitabh and Mianjy, Poorya and Mukherjee, Anirbit},
  title     = {Understanding Deep Neural Networks with Rectified Linear Units},
  booktitle = {International Conference on Learning Representations},
  year      = {2018}
}

@article{boob2022complexity,
  author  = {Boob, Digvijay and Dey, Santanu S. and Lan, Guanghui},
  title   = {Complexity of Training {ReLU} Neural Network},
  journal = {Discrete Optimization},
  volume  = {44},
  pages   = {100620},
  year    = {2022}
}

@article{dey2020approximation,
  author  = {Dey, Santanu S. and Wang, Guanyi and Xie, Yao},
  title   = {Approximation Algorithms for Training One-Node {ReLU} Neural Networks},
  journal = {IEEE Transactions on Signal Processing},
  volume  = {68},
  year    = {2020}
}

@inproceedings{froese2023training,
  author    = {Froese, Vincent and Hertrich, Christoph},
  title     = {Training Neural Networks is {NP}-Hard in Fixed Dimension},
  booktitle = {Advances in Neural Information Processing Systems},
  volume    = {36},
  year      = {2023}
}

@misc{narula2026changepointdetectionpiecewiselinear,
      title={Change-Point Detection via Piecewise Linear Fitting Using {MIP}}, 
      author={Apoorva Narula and Santanu S. Dey and Yao Xie},
      year={2026},
      eprint={2602.11947},
      archivePrefix={arXiv},
      primaryClass={math.OC},
      note = {\url{https://arxiv.org/abs/2602.11947}}, 
}

@article{goldberg2021,
  title={{MINLP} formulations for continuous piecewise linear function fitting},
  author={Goldberg, Noam and Rebennack, Steffen and Kim, Youngdae and Krasko, Vitaliy and Leyffer, Sven},
  journal={Computational Optimization and Applications},
  volume={79},
  number={1},
  pages={223--233},
  year={2021},
  publisher={Springer}
}

@article{rebennack2020,
  title={Piecewise linear function fitting via mixed-integer linear programming},
  author={Rebennack, Steffen and Krasko, Vitaliy},
  journal={INFORMS Journal on Computing},
  volume={32},
  number={2},
  pages={507--530},
  year={2020},
  publisher={INFORMS}
}

@article{cao2018multi,
  title={Multi-sensor slope change detection},
  author={Cao, Yang and Xie, Yao and Gebraeel, Nagi},
  journal={Annals of Operations Research},
  volume={263},
  number={1},
  pages={163--189},
  year={2018},
  publisher={Springer}
}

@article{fearnhead2019detecting,
  title={Detecting changes in slope with an {$\ell_0$} penalty},
  author={Fearnhead, Paul and Maidstone, Robert and Letchford, Adam},
  journal={Journal of Computational and Graphical Statistics},
  volume={28},
  number={2},
  pages={265--275},
  year={2019},
  publisher={Taylor \& Francis}
}

@article{fearnhead2024cpop,
  title={{CPOP}: Detecting changes in piecewise-linear signals},
  author={Fearnhead, Paul and Grose, Daniel},
  journal={Journal of Statistical Software},
  volume={109},
  pages={1--30},
  year={2024}
}

@book{BrentZimmermann2010,
  author    = {Brent, Richard P. and Zimmermann, Paul},
  title     = {Modern Computer Arithmetic},
  series    = {Cambridge Monographs on Applied and Computational Mathematics},
  volume    = {18},
  publisher = {Cambridge University Press},
  address   = {Cambridge},
  year      = {2010}
}

@inproceedings{wang2023holoassist,
  title={{HoloAssist}: An egocentric human interaction dataset for interactive {AI} assistants in the real world},
  author={Wang, Xin and Kwon, Taein and Rad, Mahdi and Pan, Bowen and Chakraborty, Ishani and Andrist, Sean and Bohus, Dan and Feniello, Ashley and Tekin, Bugra and Frujeri, Felipe Vieira and others},
  booktitle={2023 IEEE/CVF International Conference on Computer Vision (ICCV)},
  pages={20213--20224},
  year={2023},
  organization={IEEE}
}

@article{sagraloff2016computing,
  title={Computing real roots of real polynomials},
  author={Sagraloff, Michael and Mehlhorn, Kurt},
  journal={Journal of Symbolic Computation},
  volume={73},
  pages={46--86},
  year={2016},
  publisher={Elsevier}
}

@article{kim2011estimating,
  title={Estimating a common deterministic time trend break in large panels with cross sectional dependence},
  author={Kim, Dukpa},
  journal={Journal of Econometrics},
  volume={164},
  number={2},
  pages={310--330},
  year={2011},
  publisher={Elsevier}
}

@article{kim2014common,
  title={Common breaks in time trends for large panel data with a factor structure},
  author={Kim, Dukpa},
  journal={The Econometrics Journal},
  volume={17},
  number={3},
  pages={301--337},
  year={2014},
  publisher={Oxford University Press Oxford, UK}
}

@article{de2010breakpoints,
  title={Breakpoints in annual rainfall trends in {C{\'o}rdoba}, {Argentina}},
  author={de la Casa, Antonio and Nasello, Olga},
  journal={Atmospheric Research},
  volume={95},
  number={4},
  pages={419--427},
  year={2010},
  publisher={Elsevier}
}

@article{brilleman2017bayesian,
  title={Bayesian piecewise linear mixed models with a random change point: an application to {BMI} rebound in childhood},
  author={Brilleman, Samuel L and Howe, Laura D and Wolfe, Rory and Tilling, Kate},
  journal={Epidemiology (Cambridge, Mass.)},
  volume={28},
  number={6},
  pages={827},
  year={2017}
}

@article{wang2018high,
  title={High dimensional change point estimation via sparse projection},
  author={Wang, Tengyao and Samworth, Richard J},
  journal={Journal of the Royal Statistical Society Series B: Statistical Methodology},
  volume={80},
  number={1},
  pages={57--83},
  year={2018},
  publisher={Oxford University Press}
}

@article{enikeeva2019high,
author = {Farida Enikeeva and Zaid Harchaoui},
title = {{High-dimensional change-point detection under sparse alternatives}},
volume = {47},
journal = {The Annals of Statistics},
number = {4},
publisher = {Institute of Mathematical Statistics},
pages = {2051--2079},
year = {2019},
doi = {10.1214/18-AOS1740},
URL = {https://doi.org/10.1214/18-AOS1740}
}

@article{A2025Change,
  author = {Aghaei A., Fatemeh and Phillips, Ewan T. and Kantz, Holger},
  title = {Change point detection in {ERA5} ground temperature time series},
  journal = {Climatic Change},
  year = {2025},
  volume = {178},
  number = {11},
  pages = {193},
  doi = {10.1007/s10584-025-04021-3},
  url = {https://doi.org/10.1007/s10584-025-04021-3},
  isbn = {1573-1480}
}

@article{rahmstorf2017global,
  title={Global temperature evolution: recent trends and some pitfalls},
  author={Rahmstorf, Stefan and Foster, Grant and Cahill, Niamh},
  journal={Environmental Research Letters},
  volume={12},
  number={5},
  pages={054001},
  year={2017},
  publisher={IOP Publishing}
}

@article{beaulieu2024recent,
  title={A recent surge in global warming is not detectable yet},
  author={Beaulieu, Claudie and Gallagher, Colin and Killick, Rebecca and Lund, Robert and Shi, Xueheng},
  journal={Communications Earth \& Environment},
  volume={5},
  number={1},
  pages={576},
  year={2024},
  publisher={Nature Publishing Group UK London}
}

@article{triacca2026recent,
  title={Recent Structural Breaks in Global Temperature Series: Evidence from a Changepoint Analysis},
  author={Triacca, Umberto and Pasini, Antonello},
  journal={Climate},
  volume={14},
  number={6},
  pages={129},
  year={2026},
  publisher={MDPI}
}

@article{cahill2015change,
  title={Change points of global temperature},
  author={Cahill, Niamh and Rahmstorf, Stefan and Parnell, Andrew C},
  journal={Environmental Research Letters},
  volume={10},
  number={8},
  pages={084002},
  year={2015},
  publisher={IOP Publishing}
}

@article{foster2026global,
  title={Global warming has accelerated significantly},
  author={Foster, Grant and Rahmstorf, Stefan},
  journal={Geophysical Research Letters},
  volume={53},
  number={5},
  pages={e2025GL118804},
  year={2026},
  publisher={Wiley Online Library}
}

@book{schrijver1998theory,
  title={Theory of linear and integer programming},
  author={Schrijver, Alexander},
  year={1998},
  publisher={John Wiley \& Sons}
}

@article{weinmann2015iterative,
  title={Iterative {Potts} and {Blake--Zisserman} minimization for the recovery of functions with discontinuities from indirect measurements},
  author={Weinmann, Andreas and Storath, Martin},
  journal={Proceedings of the Royal Society A: Mathematical, Physical and Engineering Sciences},
  volume={471},
  number={2176},
  pages={20140638},
  year={2015}
}

@article{hersbach2020era5,
  author  = {Hersbach, Hans and Bell, Bill and Berrisford, Paul and Hirahara, Shoji and Hor{\'a}nyi, Andr{\'a}s and Mu{\~n}oz-Sabater, Joaqu{\'i}n and Nicolas, Julien and Peubey, Carole and Radu, Raluca and Schepers, Dinand and others},
  title   = {The {ERA5} global reanalysis},
  journal = {Quarterly Journal of the Royal Meteorological Society},
  volume  = {146},
  number  = {730},
  pages   = {1999--2049},
  year    = {2020},
  doi     = {10.1002/qj.3803}
}

@inproceedings{saxena2008turbofan,
  author    = {Saxena, Abhinav and Goebel, Kai and Simon, Don and Eklund, Neil},
  title     = {Damage propagation modeling for aircraft engine run-to-failure simulation},
  booktitle = {2008 International Conference on Prognostics and Health Management},
  pages     = {1--9},
  year      = {2008},
  organization = {IEEE}
}

@misc{boveroux2026,
  author = {Laurie Boveroux and Margarida Carvalho
            and Andrea Lodi and Quentin Louveaux},
  title = {On the Complexity of Linear Programs with
           Parametric Constraint Matrices},
  year = {2026},
  eprint = {2610.20422},
  archivePrefix = {arXiv},
  primaryClass = {math.OC},
  note = {\url{https://arxiv.org/abs/2610.20422}}
}

@inproceedings{abrahamsen2021training,
  author    = {Abrahamsen, Mikkel and Kleist, Linda and Miltzow, Tillmann},
  title     = {Training Neural Networks is {$\exists\mathbb{R}$}-complete},
  booktitle = {Advances in Neural Information Processing Systems},
  volume    = {34},
  year      = {2021}
}

@inproceedings{bertschinger2023training,
  author    = {Bertschinger, Daniel and Hertrich, Christoph and Jungeblut, Paul and Miltzow, Tillmann and Weber, Simon},
  title     = {Training Fully Connected Neural Networks is {$\exists\mathbb{R}$}-Complete},
  booktitle = {Advances in Neural Information Processing Systems},
  volume    = {36},
  pages     = {36222--36237},
  year      = {2023}
}

@article{hinkley1969inference,
  author  = {Hinkley, David V.},
  title   = {Inference about the intersection in two-phase regression},
  journal = {Biometrika},
  volume  = {56},
  number  = {3},
  pages   = {495--504},
  year    = {1969}
}

@article{feder1975asymptotic,
  author  = {Feder, Paul I.},
  title   = {On Asymptotic Distribution Theory in Segmented Regression Problems--Identified Case},
  journal = {The Annals of Statistics},
  volume  = {3},
  number  = {1},
  pages   = {49--83},
  year    = {1975}
}

@article{lerman1980fitting,
  author  = {Lerman, P. M.},
  title   = {Fitting Segmented Regression Models by Grid Search},
  journal = {Journal of the Royal Statistical Society Series C (Applied Statistics)},
  volume  = {29},
  number  = {1},
  pages   = {77--84},
  year    = {1980}
}

@article{siegmund1995using,
  author  = {Siegmund, David and Venkatraman, E. S.},
  title   = {Using the Generalized Likelihood Ratio Statistic for Sequential Detection of a Change-Point},
  journal = {The Annals of Statistics},
  volume  = {23},
  number  = {1},
  pages   = {255--271},
  year    = {1995}
}

@article{freund1985postoptimal,
  author  = {Freund, Robert M.},
  title   = {Postoptimal analysis of a linear program under simultaneous changes in matrix coefficients},
  journal = {Mathematical Programming Study},
  volume  = {24},
  pages   = {1--13},
  year    = {1985}
}

@inproceedings{lam2015numba,
  author    = {Lam, Siu Kwan and Pitrou, Antoine and Seibert, Stanley},
  title     = {Numba: A {LLVM}-based {Python} {JIT} Compiler},
  booktitle = {Proceedings of the Second Workshop on the LLVM Compiler Infrastructure in HPC},
  pages     = {1--6},
  year      = {2015},
  doi       = {10.1145/2833157.2833162}
}

@misc{gurobi,
  author = {{Gurobi Optimization, LLC}},
  title  = {Gurobi Optimizer Reference Manual},
  note   = {Version 13.0},
  year   = {2026},
  url    = {https://www.gurobi.com}
}

@article{bai2003computation,
  author  = {Bai, Jushan and Perron, Pierre},
  title   = {Computation and Analysis of Multiple Structural Change Models},
  journal = {Journal of Applied Econometrics},
  volume  = {18},
  number  = {1},
  pages   = {1--22},
  year    = {2003},
  doi     = {10.1002/jae.659}
}

@article{ROSS1976341,
title = {The arbitrage theory of capital asset pricing},
journal = {Journal of Economic Theory},
volume = {13},
number = {3},
pages = {341-360},
year = {1976},
issn = {0022-0531},
author = {Stephen A Ross},
doi = {https://doi.org/10.1016/0022-0531(76)90046-6},
url = {https://www.sciencedirect.com/science/article/pii/0022053176900466}


}

@article{bollerslev1986generalized,
  title={Generalized autoregressive conditional heteroskedasticity},
  author={Bollerslev, Tim},
  journal={Journal of econometrics},
  volume={31},
  number={3},
  pages={307--327},
  year={1986},
  publisher={Elsevier}
}

@inproceedings{guibas1978dichromatic,
  title={A dichromatic framework for balanced trees},
  author={Guibas, Leo J and Sedgewick, Robert},
  booktitle={19th Annual Symposium on Foundations of Computer Science (sfcs 1978)},
  pages={8--21},
  year={1978},
  organization={IEEE}
}

@article{hamilton1989new,
  title={A new approach to the economic analysis of nonstationary time series and the business cycle},
  author={Hamilton, James D},
  journal={Econometrica: Journal of the econometric society},
  pages={357--384},
  year={1989},
  publisher={JSTOR}
}

@article{hahn2020bayesproject,
  title={BayesProject: Fast computation of a projection direction for multivariate changepoint detection},
  author={Hahn, Georg and Fearnhead, Paul and Eckley, Idris A},
  journal={Statistics and Computing},
  volume={30},
  number={6},
  pages={1691--1705},
  year={2020},
  publisher={Springer}
}

@article{anastasiou2023generalized,
  title={Generalized multiple change-point detection in the structure of multivariate, possibly high-dimensional, data sequences},
  author={Anastasiou, Andreas and Papanastasiou, Angelos},
  journal={Statistics and Computing},
  volume={33},
  number={5},
  pages={94},
  year={2023},
  publisher={Springer}
}

%





\newpage

\appendix
\section*{Appendix}
\numberwithin{equation}{section}

{
\section{Synthetic Data Generation}
\label{app:synthetic_data}

Our synthetic data generation procedure is motivated by
multifactor models of asset returns \cite{ROSS1976341},
models of time-varying volatility
\cite{bollerslev1986generalized}, and models allowing
changes in the underlying return-generating process
\cite{hamilton1989new}. We generate $J$ correlated
stock-price trajectories, each observed over $T=100$
consecutive business days, with observation times
$x_i=i-1$ for $i\in[T]$.

\paragraph{Generating stock returns.}
Each stock is assigned independently to one of $20$
sectors and an initial price drawn uniformly from
$[20,250]$. Its daily log-return is generated as
\begin{equation}
\widetilde R_{j,t}
=
\mu_{j,t}
+\gamma_j^M R_t^M
+\gamma_j^G R_{g_j,t}^G
+\sigma_j v_{j,t}\epsilon_{j,t},
\end{equation}
where $R_t^M$ represents market-wide fluctuations,
$R_{g_j,t}^G$ represents fluctuations shared by
stocks in sector $g_j$, and $\epsilon_{j,t}$
represents independent stock-specific Gaussian
noise. The coefficients $\gamma_j^M,\gamma_j^G$
and $\sigma_j$ control the contribution of these
components. The market and stock-specific volatility
factors vary over time based on preceding shocks,
allowing periods of higher and lower variability.
Thus, stocks are correlated through shared market
and sector fluctuations while retaining
stock-specific variability.

\paragraph{Introducing changes.}
Each stock independently undergoes a change in its
baseline drift (the average daily log-return)
with probability $0.70$. For each stock, a potential
change time is sampled from a normal distribution
centered at $52\%$ of the observation period,
with a standard deviation of $10\%$ of its length.
The sampled time is rounded to an integer and
restricted to lie away from the endpoints.
The daily drift is then given by
\begin{equation}
\mu_{j,t}
=
\mu_j+\kappa_j\mathbb{I}\{t>\tau_j\},
\end{equation}
where $\tau_j$ is the change time and
$\kappa_j\sim\mathcal{N}(0,0.0015^2)$ for stocks
undergoing a change, and $\kappa_j=0$ otherwise.
Thus, different stocks may experience changes
at different times and of different magnitudes.

\paragraph{Constructing prices.}
The generated daily log-returns are capped to
$[-0.18,0.18]$ to limit extreme daily fluctuations.
Prices are then computed recursively as
\begin{equation}
P_{j,t+1}
=
P_{j,t}\exp(R_{j,t}),
\qquad
R_{j,t}
=
\operatorname{clip}(\widetilde R_{j,t},-0.18,0.18),
\end{equation}
where $\operatorname{clip}(a,l,u)
:=\min\{\max\{a,l\},u\}$.
The resulting observations are $y_i^j=P_{j,i}$.
The dense and sparse experiments use the same
generation procedure, with $J$ and $k$ specified
in Table~\ref{tab:dataset_design}.
The complete implementation, including the
parameter distributions and volatility updates,
is available at:
{\sloppy
\url{https://github.com/ApoorvaGOR/High-Dimensional-Piecewise-Linear-Fitting/blob/main/data_simulator.py}
\par}.

\paragraph{Relation to the optimization model.}
These data are not generated from the continuous
piecewise-linear model \eqref{eq:slope-change}.
Changes occur in the drift of log-returns rather
than directly in price slopes, and their locations
may differ across stocks. This allows us to
evaluate our algorithms on noisy, correlated
time series that do not necessarily satisfy
the assumed shared-breakpoint model.
}

\section{A Mixed-Integer Nonlinear Formulation and Its Big-$M$ Constants}
\label{app:minlp_bigM}

We state a mixed-integer nonlinear programming (MINLP) formulation of \eqref{eq:slope-change} for $p=2$ and describe its parameter bounds and big-$M$ constants. The construction follows \cite{narula2026changepointdetectionpiecewiselinear}.

\subsection{Formulation}

For each dimension $j\in[J]$ and segment $\ell\in\{1,2\}$, let $m_\ell^j$ and $c_\ell^j$ be the slope and the intercept of the affine function on segment $\ell$, where segment $1$ lies before the breakpoint and segment $2$ after it. In the notation of \eqref{eq:residual}, $m_1^j=\beta_j$, $c_1^j=\alpha_j$, $m_2^j=\beta_j+\delta_j$, and $c_2^j=\alpha_j-\delta_jz$. Let $\widehat y_i^j$ be the fitted value of dimension $j$ at time $x_i$, and let the binary variable $X_i$ equal $1$ if observation $i$ is assigned to segment $1$ and $0$ if it is assigned to segment $2$. The dense formulation is
\begin{align}
\min\quad&\sum_{i\in[T]}\sum_{j\in[J]}\big(y_i^j-\widehat y_i^j\big)^2\label{eq:minlp_obj}\\
\text{s.t.}\quad&-M_{1,ij}(1-X_i)\leq\widehat y_i^j-m_1^jx_i-c_1^j\leq M_{1,ij}(1-X_i),&&i\in[T],\ j\in[J],\label{eq:minlp_assign1}\\
&-M_{1,ij}X_i\leq\widehat y_i^j-m_2^jx_i-c_2^j\leq M_{1,ij}X_i,&&i\in[T],\ j\in[J],\label{eq:minlp_assign2}\\
&m_1^jz+c_1^j=m_2^jz+c_2^j,&&j\in[J],\label{eq:minlp_continuity}\\
&x_i\leq z+M_{2,i}(1-X_i),&&i\in[T],\label{eq:minlp_localization_first}\\
&x_i\geq z-M_{3,i}X_i,&&i\in[T],\label{eq:minlp_localization_second}\\
&X_i\geq X_{i+1},&&i\in[T-1],\label{eq:minlp_monotonicity}\\
&(m,c,z)\in\Theta,\qquad X\in\{0,1\}^T,\label{eq:minlp_domain}
\end{align}
where the set $\Theta$ and the constants $M_{1,ij}$, $M_{2,i}$, and $M_{3,i}$ are defined below. The constraints \eqref{eq:minlp_assign1} and \eqref{eq:minlp_assign2} set $\widehat y_i^j$ to the value of the affine function of the segment to which observation $i$ is assigned. The constraints \eqref{eq:minlp_localization_first}--\eqref{eq:minlp_monotonicity} ensure that the observations of segment $1$ are exactly those up to the breakpoint, with either assignment allowed for an observation at the breakpoint. The constraints \eqref{eq:minlp_continuity} impose continuity at the breakpoint. They are bilinear, and they are the source of nonconvexity. For a fixed breakpoint $z$, every feasible point corresponds to a continuous piecewise-linear fit with breakpoint $z$ in each dimension, and the objective \eqref{eq:minlp_obj} equals $L_2(r)$.

In the sparse problem, we add a binary variable $\xi_j$ for each $j\in[J]$, where $\xi_j=1$ indicates that dimension $j$ may change slope, and the constraints of Section~\ref{app:minlp_sparse}.

\subsection{Parameter bounds}

For each dimension $j$, the lower and upper slope bounds are the smallest and largest slopes between two observations,
\begin{align}
L_m^j&:=\min_{\substack{i,i'\in[T]\\i\neq i'}}\frac{y_i^j-y_{i'}^j}{x_i-x_{i'}},\label{eq:minlp_Lm}\\
U_m^j&:=\max_{\substack{i,i'\in[T]\\i\neq i'}}\frac{y_i^j-y_{i'}^j}{x_i-x_{i'}}.\label{eq:minlp_Um}
\end{align}
Given these slope bounds, the intercept bounds are
\begin{align}
L_c^j&:=\min_{i\in[T]}\min\big\{y_i^j-L_m^jx_i,\ y_i^j-U_m^jx_i\big\},\label{eq:minlp_Lc}\\
U_c^j&:=\max_{i\in[T]}\max\big\{y_i^j-L_m^jx_i,\ y_i^j-U_m^jx_i\big\}.\label{eq:minlp_Uc}
\end{align}
For every $j\in[J]$ and $\ell\in\{1,2\}$, we impose
\begin{equation}
L_m^j\leq m_\ell^j\leq U_m^j,\qquad L_c^j\leq c_\ell^j\leq U_c^j.\label{eq:minlp_parameter_bounds}
\end{equation}
The breakpoint is restricted to $x_2\leq z\leq x_{T-1}$, as in \eqref{eq:slope-change}, so that at least two observations lie on each side of the breakpoint. Hence
\begin{equation}
\Theta:=\Big\{(m,c,z):\ L_m^j\leq m_\ell^j\leq U_m^j,\ L_c^j\leq c_\ell^j\leq U_c^j\ \text{for all } j\in[J],\ \ell\in\{1,2\},\ x_2\leq z\leq x_{T-1}\Big\}.\label{eq:minlp_theta}
\end{equation}

The bounds defining $\Theta$ are data-driven, and they are not guaranteed to contain an optimal solution of \eqref{eq:slope-change}, whose slopes and intercepts are unrestricted. Hence the formulation \eqref{eq:minlp_obj}--\eqref{eq:minlp_domain} is a restriction of \eqref{eq:slope-change}, and its optimal value is at least $\textup{OPT}_{2,k}$. The two values are equal if some optimal solution of \eqref{eq:slope-change} lies in $\Theta$. This can be checked using the solution returned by the algorithms of Sections~\ref{sec:nonsparse} and~\ref{sec:sparse}.

\subsection{Big-$M$ constants for the fitted values}

For $i\in[T]$ and $j\in[J]$, let
\begin{align}
\underline f_{ij}&:=\min\big\{L_m^jx_i,\ U_m^jx_i\big\}+L_c^j,\label{eq:minlp_f_lower}\\
\overline f_{ij}&:=\max\big\{L_m^jx_i,\ U_m^jx_i\big\}+U_c^j.\label{eq:minlp_f_upper}
\end{align}
For every $(m,c,z)\in\Theta$ and $\ell\in\{1,2\}$, we have $m_\ell^jx_i\in\big[\min\{L_m^jx_i,U_m^jx_i\},\ \max\{L_m^jx_i,U_m^jx_i\}\big]$, and therefore
\begin{equation}
\underline f_{ij}\leq m_\ell^jx_i+c_\ell^j\leq\overline f_{ij}.\label{eq:minlp_affine_range}
\end{equation}
Let
\begin{equation}
M_{1,ij}:=\overline f_{ij}-\underline f_{ij}.\label{eq:minlp_M1}
\end{equation}
At every feasible point, $\widehat y_i^j$ equals $m_\ell^jx_i+c_\ell^j$ for the segment $\ell$ to which observation $i$ is assigned, so $\underline f_{ij}\leq\widehat y_i^j\leq\overline f_{ij}$ by \eqref{eq:minlp_affine_range}. Hence, for the other segment $\ell'$, \eqref{eq:minlp_affine_range} gives
\begin{equation}
\big|m_{\ell'}^jx_i+c_{\ell'}^j-\widehat y_i^j\big|\leq M_{1,ij}.\label{eq:minlp_M1_validity}
\end{equation}
Thus the assignment constraint of the segment to which observation $i$ is not assigned is satisfied automatically, and $M_{1,ij}$ is a valid big-$M$ constant in \eqref{eq:minlp_assign1} and \eqref{eq:minlp_assign2}.

\subsection{Big-$M$ constants for the breakpoint}

Since $x_2\leq z\leq x_{T-1}$, valid constants for \eqref{eq:minlp_localization_first} and \eqref{eq:minlp_localization_second} are
\begin{align}
M_{2,i}&:=\max\{0,x_i-x_2\},&&i\in[T],\label{eq:minlp_M2}\\
M_{3,i}&:=\max\{0,x_{T-1}-x_i\},&&i\in[T].\label{eq:minlp_M3}
\end{align}
If $X_i=1$, then \eqref{eq:minlp_localization_first} enforces $x_i\leq z$, and \eqref{eq:minlp_localization_second} reads $x_i\geq z-M_{3,i}$, which holds since $z\leq x_{T-1}$. If $X_i=0$, then \eqref{eq:minlp_localization_second} enforces $x_i\geq z$, and \eqref{eq:minlp_localization_first} reads $x_i\leq z+M_{2,i}$, which holds since $z\geq x_2$. The monotonicity constraints \eqref{eq:minlp_monotonicity} ensure that the observations of each segment form a contiguous block.

\subsection{Constraints and Big-$M$ constants for the sparse problem}
\label{app:minlp_sparse}

The bounds \eqref{eq:minlp_parameter_bounds} imply
\begin{equation}
|m_2^j-m_1^j|\leq U_m^j-L_m^j,\qquad |c_2^j-c_1^j|\leq U_c^j-L_c^j.\label{eq:minlp_differences}
\end{equation}
Hence valid big-$M$ constants are
\begin{equation}
M_m^j:=U_m^j-L_m^j,\qquad M_c^j:=U_c^j-L_c^j,\label{eq:minlp_sparse_M}
\end{equation}
and the sparse constraints are
\begin{align}
-M_m^j\xi_j&\leq m_2^j-m_1^j\leq M_m^j\xi_j,&&j\in[J],\label{eq:minlp_sparse_slopes}\\
-M_c^j\xi_j&\leq c_2^j-c_1^j\leq M_c^j\xi_j,&&j\in[J],\label{eq:minlp_sparse_intercepts}\\
\sum_{j\in[J]}\xi_j&\leq k,\label{eq:minlp_sparse_cardinality}\\
\xi_j&\in\{0,1\},&&j\in[J].\label{eq:minlp_sparse_binary}
\end{align}
If $\xi_j=0$, these constraints enforce $m_1^j=m_2^j$ and $c_1^j=c_2^j$, so dimension $j$ is fitted by a single line, that is, $\delta_j=0$. If $\xi_j=1$, they impose no restriction beyond \eqref{eq:minlp_parameter_bounds}.

\end{document}